\documentclass[11pt,a4paper]{article}

\usepackage{epsf,epsfig,amsfonts,amsgen,amsmath,amstext,amsbsy,amsopn,amsthm}
\usepackage{amsmath,times,mathptmx}
\usepackage{amsfonts,amsthm,amssymb}
\usepackage{amsfonts}
\usepackage{graphics}
\usepackage{latexsym,bm}
\usepackage{amsfonts,amsthm,amssymb,bbding}
\usepackage{indentfirst}
\usepackage{graphicx}
\usepackage{color}
\usepackage[colorlinks=true,anchorcolor=blue,filecolor=blue,linkcolor=blue,urlcolor=blue,citecolor=blue]{hyperref}
\usepackage{float}
\usepackage{tikz,enumerate}
\evensidemargin=\oddsidemargin\topmargin=-1.5cm

\newtheorem{thm}{Theorem}[section]

\newtheorem{lem}{Lemma}[section]
\newtheorem{cor}{Corollary}[section]

\newtheorem{claim}{Claim}[section]
\newtheorem{definition}{Definition}[section]
\renewcommand\proofname{\bf Proof}
\addtocounter{section}{0}

\begin{document}
\title{\bf On the Tur\'{a}n number of odd-ballooning of  $3$-chromatic graphs\footnote{Supported by the National Natural Science Foundation of China (Nos. \!12501471, \!12271162, \!12171066, \!12471311, \!12471324), the Natural Science Foundation of Shanghai (Nos. \!22ZR1416300, \!24ZR1415500),
and the Research and Innovation Team of Anhui Province (No. 2024AH010023).}}
\author{ {\bf Longfei Fang$^{a,b}$},~{\bf Xueyi Huang$^{a}$}\thanks{Corresponding author.} \setcounter{footnote}{-1}\footnote{\emph{Email address:} lffang@chzu.edu.cn (L. Fang), huangxy@ecust.edu.cn (X. Huang), huiqiulin@126.com (H. Lin), jlshu@shnu.edu.cn (J. Shu).},~
{\bf Huiqiu Lin$^{a}$},~
{\bf Jinlong Shu$^{c}$}
\\[4mm]
\small $^{a}$School of Mathematics, East China University of Science and Technology,\\
\small  Shanghai 200237, China\\
\small $^{b}$School of Mathematics and Finance, Chuzhou University, \\
\small  Chuzhou, Anhui 239012, China\\
\small $^{c}$School of Finance Business, Shanghai Normal University,\\
\small   Shanghai 200234, China\\
}

\date{}
\maketitle
{\flushleft\large\bf Abstract}
Given a graph $F$, the Tur\'{a}n number ${\rm ex}(n,F)$ is the maximum number of edges
in any $n$-vertex $F$-free graph. The odd-ballooning of $F$, denoted by $F^{o}$,
is a graph obtained by replacing each edge of $F$ with an odd cycle, where all new vertices
of the odd cycles are distinct. The Tur\'{a}n number of the odd-ballooning of $F$ has been established for several
important cases. For a star, it was determined by Erd\H{o}s, F\"{u}redi, Gould, and
Gunderson (1995), Hou, Qiu, and Liu (2018), and Yuan (2018);
for trees under certain conditions, by Zhu and Chen (2023); and for complete
bipartite graphs $K_{s,t}$ ($t\geq s \geq 2$) where each substituted odd cycle has length at least five,
by Peng and Xia (2024). In this paper, we apply Simonovits' celebrated method of progressive induction to determine
the Tur\'{a}n number for the odd-ballooning of a class of $3$-chromatic graphs. Specifically,
let $F$ be a graph formed by connecting a single vertex to all vertices of another graph
whose components are either non-trivial trees or even cycles. We determine ${\rm ex}(n,F^{o})$
when each substituted odd cycle in $F^{o}$ has length at least five. As corollaries, we obtain the Tur\'{a}n number for the odd-ballooning of several well-known graph classes, including odd wheels, fan graphs, book graphs, and friendship graphs, where each substituted odd cycle in the ballooning has length at least five.

\begin{flushleft}
\textbf{Keywords:} Tur\'{a}n number; odd-ballooning; progressive induction; decomposition family

\end{flushleft}
\textbf{AMS Classification:} 05C35

\section{Introduction}

Let $\mathcal{F}$ be a family of graphs.
A graph is called \emph{$\mathcal{F}$-free}
if it contains no subgraph isomorphic to any member $F\in \mathcal{F}$.
The Tur\'{a}n number of $\mathcal{F}$, denoted by ${\rm ex}(n,\mathcal{F})$,
is the maximum number of edges in any $n$-vertex $\mathcal{F}$-free graph.
Let ${\rm EX}(n,\mathcal{F})$ denote the family of $\mathcal{F}$-free graphs
with  $n$ vertices and ${\rm ex}(n,\mathcal{F})$ edges.
When $\mathcal{F}$ consists of a single graph $F$,
we use ${\rm ex}(n,F)$ (resp. ${\rm EX}(n,F)$) instead of ${\rm ex}(n,\mathcal{F})$ (resp. ${\rm EX}(n,\mathcal{F})$). As one of the earliest results in extremal graph theory,
Tur\'{a}n's theorem \cite{Turan} states that
${\rm EX}(n,K_{r+1})=\{T_{n,r}\}$,
where the Tur\'{a}n graph $T_{n,r}$ is the complete $r$-partite graph on $n$ vertices
with part sizes as equal as possible.
For the Tur\'{a}n number of general graphs,
the celebrated  Erd\H{o}s-Stone-Simonovits theorem \cite{Erdos-1966,ES1946} states that
$${\rm ex}(n,F)=\Big(1-\frac{1}{\chi(F)-1}\Big)\binom{n}{2}+o(n^2),$$
where $\chi(F)$ denotes the chromatic number of $F$.
When $\chi(F) \geq 3$, this theorem gives the asymptotic value of $\mathrm{ex}(n,F)$.
Determining the exact value of $\mathrm{ex}(n,F)$, however, remains an important and challenging open problem in extremal graph theory.

For a simple graph $G$, we write $V(G)$ for the vertex set, $E(G)$ for the edge set, $|G| := |V(G)|$ for the number of vertices, and $e(G) := |E(G)|$ for the number of edges. For a positive integer $k$, let $K_{k}$, $E_k$, $S_k$, and $C_k$ denote the complete graph,
the empty graph, the star, and the cycle on $k$ vertices, respectively. Let $F \nabla H$ denote the  join and $F \cup H$ the disjoint union of graphs $F$ and $H$. To study the Tur\'{a}n number of $F$, we introduce the \emph{decomposition family} $\mathcal{M}(F)$,
a concept originally developed by Simonovits. While Simonovits first presented a similar notion in \cite{Simonovits-1974},
the specific term ``decomposition family'' was not used until later in \cite{Simonovits-1983}
(see also \cite{Simonovits-2019}). The version of the decomposition family as defined here
was given by Liu \cite{LIU}.

\begin{definition}\label{def2.1}
Given a graph $F$ with $\chi(F)=r+1\geq 3$,
the decomposition family $\mathcal{M}(F)$ consists of minimal graphs $M$
such that $F\subseteq (M\cup E_t)\nabla T_{(r-1)t,r-1}$ for some constant $t=t(F)$.
\end{definition}

This definition provides a framework for characterizing extremal graphs for $F$ using the Tur\'{a}n graph $T_{n,r}$.
The approach suggests embedding a maximal $\mathcal{M}(F)$-free graph into one part of $T_{n,r}$ to construct extremal graphs for $F$.
Remarkably, this construction yields the exact Tur\'{a}n number $\mathrm{ex}(n,F)$ for many graphs $F$. Let $F_t$ denote the \emph{friendship graph} formed by $t$ triangles sharing a common vertex.
For integers $\nu$ and $\Delta$, we define the extremal function
\[
f(\nu,\Delta) = \max\{e(G) : \nu(G) \leq \nu, \Delta(G) \leq \Delta\},
\]
where $\nu(G)$ and $\Delta(G)$ represent the matching number and maximum degree of $G$, respectively.
Abbott, Hanson, and Sauer \cite{Abbott-1972} determined the exact value of $f(t-1,t-1)$. The following theorem by Erd\H{o}s, F\"{u}redi, Gould, and Gunderson \cite{EFGG1995} serves as a classical illustration of this approach.

\begin{thm}\label{thm1.2}\emph{(\cite{EFGG1995})}
For each $t\geq 1$ and each $n\geq 50t^2$,
${\rm ex}(n,F_t)=e(T_{n,2})+f(t-1,t-1)$.
\end{thm}

A matching $M_t$ is the vertex-disjoint union of $t$ edges.
Note that $\mathcal{M}(F_t)=\{M_t,S_{t+1}\}$ and ${\rm ex}(n,\{M_t,S_{t+1}\})=f(t-1,t-1)$.
Therefore, by embedding a maximal $\{M_t,S_{t+1}\}$-free graph into one part of  $T_{n,2}$,
we obtain an extremal graph for $F_t$.

Theorem \ref{thm1.2} has been extended in several directions.
Chen, Gould, Pfender, and Wei \cite{Chen-2003} proved that
\[
\mathrm{ex}(n,K_1 \nabla tK_r) = e(T_{n,r}) + f(t-1,t-1)
\]
for all integers $t \geq 1$, $r \geq 2$, and $n \geq 16t^3(r+1)^8$. Building on these results, Liu \cite{LIU} introduced the notion of edge blow-up of graphs.
For a graph $F$ and an integer $r \geq 2$, the \emph{edge blow-up} $F^{r+1}$ is obtained by replacing each edge of $F$ with a copy of $K_{r+1}$,
where all new vertices are distinct.
Note that $K_1 \nabla tK_r \cong S_{t+1}^{r+1}$. To date, the exact Tur\'{a}n number $\mathrm{ex}(n,F^{r+1})$ has been determined for numerous graphs.
For more information on this topic,
we refer the readers to \cite{E1962,LIU,Moon1968,Simonovits1968,Song2024,Wang2021,Y2022}.

If we consider $K_3$ as an odd cycle, then Theorem \ref{thm1.2} admits another natural generalization.
The \emph{odd-ballooning} of a graph $F$, denoted by $F^{o}$, is obtained by replacing each edge of $F$ with an odd cycle, where all new vertices from different substituted cycles are distinct.
Note that the friendship graph $F_t$ is precisely the odd-ballooning of the star $S_{t+1}$. The Tur\'{a}n number for the odd-ballooning of graphs was first investigated by Hou, Qiu, and Liu \cite{Hou2016,Hou2018} and Yuan \cite{Yuan2018} for the case of stars.
Zhu, Kang, and Shan \cite{Zhu2020} subsequently determined the Tur\'{a}n number for the
odd-ballooning of paths and cycles.
Further progress was made by Zhu and Chen \cite{Zhu2023}, who obtained the Tur\'{a}n number for the odd-ballooning of trees under certain conditions, generalizing previous results for paths and stars.
Most recently, Peng and Xia \cite{Peng-2024}  determined the Tur\'{a}n number for the odd-ballooning of complete bipartite graphs where each substituted cycle has length at least five.

Simonovits' theorem, as presented in \cite[Theorem 1]{Simonovits-1974}, asserts that when the decomposition family of the forbidden graph $F$ includes a linear forest,
the extremal graphs for $F$ have simple and symmetric structures. Due to the theorem's substantial complexity, we refer interested readers to \cite{Simonovits-1974} for complete details. We remark that the decomposition family of the odd-ballooning of the aforementioned graphs always contains a matching.

{\flushleft\textbf{Notation.}}  Throughout this paper,
$F^{o}$ denotes the odd-ballooning of a graph $F$, where each edge is replaced by an odd cycle of length at least five. Note that different edges may be replaced by odd cycles of varying lengths.

We begin by presenting the following result for graphs $F$ with chromatic number $\chi(F) \geq 4$.

\begin{thm}\label{thm1.2A}
For sufficiently large $n$ and $\chi(F)\geq 4$, we have ${\rm EX}(n,F^{o})={\rm EX}(n,F)$.
\end{thm}

This result naturally leads us to investigate the Tur\'{a}n number $\mathrm{ex}(n,F^{o})$ for graphs $F$ with $\chi(F) \in \{2,3\}$.  A particularly interesting class consists of graphs $F$ formed by connecting a single vertex to all vertices of another graph whose components are either non-trivial trees or even cycles.
In contrast to previous research on $F^{o}$, where the presence of linear forests in $\mathcal{M}(F^{o})$ underpins the structural analysis, the family $\mathcal{M}(F^{o})$ considered in this paper lacks such structures. This distinction constitutes the primary novelty and difficulty of our study, making the exact evaluation of $\mathrm{ex}(n,F^{o})$ a meaningful and challenging problem. Our main result is stated as follows.

\begin{thm}\label{thm1.1}
Let $n$ be sufficiently large,
and let $F=K_1\nabla F^{\bullet}$ where each component of $F^{\bullet}$ is either a non-trivial tree or an even cycle. Then
 $${\rm EX}(n,F^{o})=
   \begin{cases}
   \{(M_{\frac{1}{2}|F^{\bullet}|-1}\cup K_1)\nabla T_{n-|F^{\bullet}|+1,2}\},
       & \hbox{if each component of $F^{\bullet}$ is an even cycle;} \\
     \{K_{e(F^{\bullet})}\nabla T_{n-e(F^{\bullet}),2}\},
       & \hbox{if each component of $F^{\bullet}$ is an edge;} \\
          \{E_{e(F^{\bullet})}\nabla T_{n-e(F^{\bullet}),2}\},
            & \hbox{otherwise.}
              \end{cases}
$$
\end{thm}

Let $W_{2k+1}=K_1\nabla C_{2k}$ denote the \emph{odd wheel},
$H_k=K_1\nabla P_{k+1}$ the \emph{fan graph},
and $B_k=K_1\nabla S_{k+1}$ the \emph{book graph}.
By Theorem \ref{thm1.1}, we immediately deduce the following corollary.

\begin{cor}\label{cor1.1}
Let $k\geq 2$ and $n$ be sufficiently large. Then
\begin{enumerate}[{\rm (i)}]\setlength{\itemsep}{0pt}
\item ${\rm EX}(n,W_{2k+1}^{o})=\{(M_{k-1}\cup K_1)\nabla T_{n-2k+1,2}\}$;
\item ${\rm EX}(n,H_{k}^{o})=\{E_{k}\nabla T_{n-k,2}\}$;
\item ${\rm EX}(n,B_{k}^{o})=\{E_k\nabla T_{n-k,2}\}$;
\item ${\rm EX}(n,F_{k}^{o})=\{K_k\nabla T_{n-k,2}\}$.
\end{enumerate}
\end{cor}

The remainder of this paper is organized as follows. In Section \ref{sec2}, we characterize the decomposition family of the odd-ballooning of graphs and present several relevant lemmas. Section \ref{sec3} contains technical lemmas needed for our analysis. Finally, we present the proofs of our main results: Theorem \ref{thm1.2A} in Section \ref{sec5} and Theorem \ref{thm1.1} in Section \ref{sec4}.

\section{Decomposition family of the odd-ballooning of graphs} \label{sec2}

To state our main theorem and related results, we first introduce some definitions and notations.
Given two disjoint vertex subsets $X,Y\subseteq V(G)$.
Let $G[X]$ be the subgraph induced by $X$, $G-X=G[V(G)\setminus X]$,
and $G[X,Y]$ be the bipartite subgraph on the vertex set $X\cup Y$
which consists of all edges with one
endpoint in $X$ and the other in $Y$.
For short, we write $e(X)=e(G[X])$ and $e(X,Y)=e(G[X,Y])$, respectively.

\begin{figure}[!h]
\centering
\begin{tikzpicture}[scale=0.75, x=1.00mm, y=1.00mm, inner xsep=0pt, inner ysep=0pt, outer xsep=0pt, outer ysep=0pt]
\definecolor{L}{rgb}{0,0,0}
\definecolor{F}{rgb}{0,0,0}

\draw [line width=0.3mm] (0,0) circle (1.5cm);
\draw(-2,-1) node[anchor=base west]{\fontsize{10.23}{17.07}\selectfont $F-\{u\}$};

\node[circle,fill=F,draw,inner sep=0pt,minimum size=1.2mm] (u) at (-30.00,0.00) {};
\draw(-34,-1) node[anchor=base west]{\fontsize{10.23}{17.07}\selectfont $u$};

\node[circle,fill=F,draw,inner sep=0pt,minimum size=1.2mm] (v1) at (-10.00,8.00) {};
\draw(-7,7) node[anchor=base west]{\fontsize{10.23}{17.07}\selectfont $v_1$};

\node[circle,fill=F,draw,inner sep=0pt,minimum size=1.2mm] (v2) at (-10.00,4.00) {};
\draw(-7,3) node[anchor=base west]{\fontsize{10.23}{17.07}\selectfont $v_2$};

\node[circle,fill=F,draw,inner sep=0pt,minimum size=0.3mm] () at (-10.00,0.00) {};
\node[circle,fill=F,draw,inner sep=0pt,minimum size=0.3mm] () at (-10.00,-2.00) {};
\node[circle,fill=F,draw,inner sep=0pt,minimum size=0.3mm] () at (-10.00,-4.00) {};

\node[circle,fill=F,draw,inner sep=0pt,minimum size=1.2mm] (vd) at (-10.00,-8.00) {};
\draw(-7,-9) node[anchor=base west]{\fontsize{10.23}{17.07}\selectfont $v_d$};

\definecolor{L}{rgb}{0,0,0}
\path[line width=0.3mm, draw=L] (u) -- (v1);
\path[line width=0.3mm, draw=L] (u) -- (v2);
\path[line width=0.3mm, draw=L] (u) -- (vd);

\draw[-latex][line width=0.4mm]  (25,0) -- (35,0);


\draw [line width=0.3mm] (80,0) circle (1.5cm);
\draw(78,-1) node[anchor=base west]{\fontsize{10.23}{17.07}\selectfont $F-\{u\}$};

\node[circle,fill=F,draw,inner sep=0pt,minimum size=1.2mm] (u1R) at (50.00,8.00) {};
\draw(42,7) node[anchor=base west]{\fontsize{10.23}{17.07}\selectfont $u_1$};
\node[circle,fill=F,draw,inner sep=0pt,minimum size=1.2mm] (u2R) at (50.00,4.00) {};
\draw(42,3) node[anchor=base west]{\fontsize{10.23}{17.07}\selectfont $u_2$};
\node[circle,fill=F,draw,inner sep=0pt,minimum size=0.3mm] () at (50.00,0.00) {};
\node[circle,fill=F,draw,inner sep=0pt,minimum size=0.3mm] () at (50.00,-2.00) {};
\node[circle,fill=F,draw,inner sep=0pt,minimum size=0.3mm] () at (50.00,-4.00) {};
\node[circle,fill=F,draw,inner sep=0pt,minimum size=1.2mm] (udR) at (50.00,-8.00) {};
\draw(42,-9) node[anchor=base west]{\fontsize{10.23}{17.07}\selectfont $u_d$};
\node[circle,fill=F,draw,inner sep=0pt,minimum size=1.2mm] (wdR) at (60.00,-8.00) {};
\draw(58,-12) node[anchor=base west]{\fontsize{10.23}{17.07}\selectfont $w_d$};
\draw(49,-18) node[anchor=base west]{\fontsize{10.23}{17.07}\selectfont $I_u$};

\node[circle,fill=F,draw,inner sep=0pt,minimum size=1.2mm] (v1R) at (70.00,8.00) {};
\draw(73,7) node[anchor=base west]{\fontsize{10.23}{17.07}\selectfont $v_1$};

\node[circle,fill=F,draw,inner sep=0pt,minimum size=1.2mm] (v2R) at (70.00,4.00) {};
\draw(73,3) node[anchor=base west]{\fontsize{10.23}{17.07}\selectfont $v_2$};

\node[circle,fill=F,draw,inner sep=0pt,minimum size=0.3mm] () at (70.00,0.00) {};
\node[circle,fill=F,draw,inner sep=0pt,minimum size=0.3mm] () at (70.00,-2.00) {};
\node[circle,fill=F,draw,inner sep=0pt,minimum size=0.3mm] () at (70.00,-4.00) {};

\node[circle,fill=F,draw,inner sep=0pt,minimum size=1.2mm] (vdR) at (70.00,-8.00) {};
\draw(73,-9) node[anchor=base west]{\fontsize{10.23}{17.07}\selectfont $v_d$};

\definecolor{L}{rgb}{1,0,0}
\path[line width=0.3mm, draw=L] (u1R) -- (v1R);
\path[line width=0.3mm, draw=L] (u2R) -- (v2R);
\definecolor{L}{rgb}{0,0,1}
\path[line width=0.3mm, draw=L] (udR) -- (wdR);

\definecolor{L}{rgb}{1,0,1}
\path[line width=0.20mm, draw=L,dashed] (50,0) ellipse (3mm and 12mm);

\end{tikzpicture}
\caption{A cracking of $F$ on the vertex $u$.}{\label{fig-1.1}}
\end{figure}
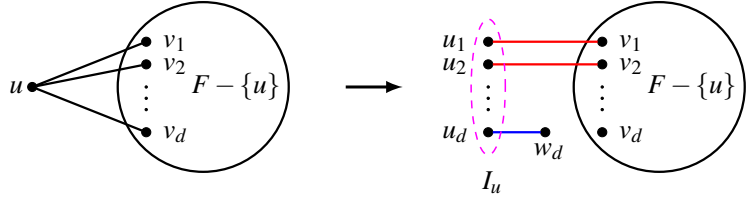

Let $F$ be a graph with $\chi(F)\in \{2,3\}$, and $u$ be a vertex of $F$ with $N_F(u)=\{v_1,\dots,v_d\}$.
A \emph{cracking} on $u$ involves replacing
$u$ with an independent set $I_u:=\{u_1,u_2,\dots,u_d\}$,
and replacing $uv_i$ with a new edge $e(u,v_i)$ for each $i\in [d]$, where $[d]:=\{1,2,\dots,d\}$.
For each $i$, $e(u,v_i)$ is defined by the following mutually exclusive cases:
\begin{itemize}\setlength{\itemsep}{0pt}
 \item Case 1: $e(u,v_i):=u_iv_i$. In this case, $uv_i$ is classified as Type I;
\item Case 2: $e(u,v_i):=u_iw_i$, where $w_i$ is a new vertex.
Here,  $uv_i$ is classified as Type II.
\end{itemize}
Additionally,  all new vertices must be distinct, see Figure \ref{fig-1.1}.
Let $U$ be an independent set of $F$.
We use $\mathcal{C}(F,U)$ to denote the family of graphs that can be obtained from $F$
by cracking all vertices in $U$.
Define $$\mathcal{C}(F):=\{\mathcal{C}(F,U)~|~U~\text{is an independent set of}~F\}.$$

For any edge $uv\in E(F)$,
let $C(u,v)$ denote the substituted odd cycle associated with $uv$ in $F^{o}$,
and define $\ell_{u,v}:=|C(u,v)|$.
Now we characterize the decomposition family of $F^{o}$.

\begin{lem}\label{lem2.1}
Let $F$ be a graph with $\chi(F)\in \{2,3$\}.
Then $\mathcal{C}(F)=\mathcal{M}(F^{o})$.
\end{lem}

\begin{proof}
Let $n$ be sufficiently large, and let $X$ and $Y$ be two color classes of $T_{n,2}$.
Take an arbitrary $J\in \mathcal{M}(F^{o})$,
and let $G$ be the graph obtained from $T_{n,2}$ by embedding $J$ into $X$.
By the definition of $\mathcal{M}(F^{o})$, $G$ contains a subgraph $F^{o}$.
Moreover, since $F$ is the skeleton of $F^{o}$,
it follows that $F$ is also a subgraph of $G$.

We first prove that $\mathcal{M}(F^{o})\subseteq \mathcal{C}(F)$.
Select an arbitrary edge $uv\in E(F)$.
Since $T_{n,2}$ contains no odd cycles,
it follows that $|E(J)\cap E(C(u,v))|\geq 1$.
Suppose that $|E(J)\cap E(C(u,v))|\geq 2$,
that is, $C(u,v)$ contains exactly $s\geq 2$ edges $u_1v_1,u_2v_2,\dots,u_sv_s$ in $E(J)$,
where $u_1v_1=uv$ if $uv\in E(J)$.
Now, we select a new cycle $C^*(u,v)$ of length $\ell_{u,v}$ such that $E(C^*(u,v))\cap E(J)=\{u_1v_1\}$ and $V(C^*(u,v))\cap V(F^{o})=\{u,v\}\cup \{u_1,v_1\}$.
Let $G'$ be the graph obtained from $T_{n,2}$ by embedding $J'$ into $X$,
where $J'=J-\{u_2v_2,\dots,u_sv_s\}$.
From the construction of $C^*(u,v)$ we know that $G'$ contains a copy of $F^{o}$.
This is impossible as $J\in \mathcal{M}(F^{o})$ and $J'$ is a proper subgraph of $J$.
Thus, $|E(J)\cap E(C(u,v))|=1$.
In other words, for every $uv\in E(F)$,
$C(u,v)$ contributes precisely one edge, denoted as $e_{u,v}$, to $E(J)$.
Thus $$|E(J)|=|E(F)|~\text{and}~E(J)=\{e_{u,v}~:~uv\in E(F)\}.$$

Recall that $F$ is a subgraph of $G$.
Set $X_1=V(F)\cap X$ and $Y_1=V(F)\cap Y$.
Then, $X_1\cup Y_1$ is a partition of $V(F)$.
Let $u$ be an arbitrary vertex in $Y_1$ and $N_{F}(u)=\{v_1,v_2,\dots,v_d\}$.
Since $Y_1$ is an independent set of $F$, we get $\{v_1,v_2,\dots,v_d\}\subseteq X_1$.
Now we select a graph $J''\in \mathcal{C}(F,Y_1)$ such that for each $uv_i\in E(F[X_1,Y_1])$,
we define $e(u,v_i)\in E(J'')$ as follows:
\begin{itemize}\setlength{\itemsep}{0pt}
\item $e(u,v_i):=u_iv_i$ if $e_{u,v_i}$ is incident to $v_i$;
\item $e(u,v_i):=u_iw_i$ if $e_{u,v_i}$ is not incident to $v_i$.
\end{itemize}
Now, let $e\mapsto e$ if $e\in E(F[X_1])$,
and $e_{u,v_i}\mapsto e(u,v_i)$ if $uv_i\in E(F[X_1,Y_1])$.
This means that $J\cong J''$.
Since $J$ is arbitrary, we conclude that $\mathcal{M}(F^{o})\subseteq \mathcal{C}(F)$.

Conversely, we shall prove that $\mathcal{C}(F)\subseteq \mathcal{M}(F^{o})$.
Let $J^*$ be an arbitrary graph in $\mathcal{C}(F)$.
Then there exists an independent set $Y_1$ of $F$ such that $J^*\in \mathcal{C}(F,Y_1)$.
Set $X_1=V(F)\setminus Y_1$.
Clearly, $X_1=V(F)\cap V(J^*)$.
Let $G^*$ be the graph obtained from $T_{n,2}$ by embedding $J^*$ into the color class $X$ and $Y_1$ into the color class $Y$ of $T_{n,2}$.
Observe that $|E(J^*)|=|E(F)|$, and $|E(J)|=|E(F)|$ for all $J\in \mathcal{M}(F^{o})$.
Consequently, to prove $J^*\in \mathcal{M}(F^{o})$,
it is sufficient to show that $G^*$ contains a copy of $F^{o}$.
Let $uv$ be any edge of $F[X_1]$.
Then $G^*$ contains a cycle $C^*(u,v)$ of length $\ell_{u,v}$ such that $E(C^*(u,v))\cap E(J^{*})=\{uv\}$.
Let $uv_i$ be any edge of $ F[X_1,Y_1]$ with $u\in Y_1$ and $v_i\in X_1$.
If $uv_i$ is of Type I,
then $G^*$ contains a cycle $C^*(u,v_i)$ of length $\ell_{u,v_i}$ such that
$E(C^*(u,v_i))\cap E(J^{*})=\{u_iv_i\}$;
If $uv_i$ is of Type II,
then $G^*$ contains a cycle $C^*(u,v_i)$ of length $\ell_{u,v_i}$ such that $E(C^*(u,v_i))\cap E(J^{*})=\{u_iw_i\}$.
In addition, by appropriately selecting vertices in all odd cycles mentioned above,
$G^*$ contains a copy of $F^{o}$.
Thus, $\mathcal{C}(F)\subseteq \mathcal{M}(F^{o})$.

Combining the two inclusions, we conclude that $\mathcal{C}(F)=\mathcal{M}(F^{o})$, as needed.
\end{proof}

A \emph{covering} of a graph is a set that meets all edges.
We use $\beta(F)$ to denote the minimum number of vertices in a covering of $F$.
An \emph{independent set} of a graph is a set of vertices where no two vertices are adjacent.
Similarly, an \emph{independent covering} of a bipartite graph is an independent set which meets all edges.
Let $\mathcal{F}$ be a finite graph family with $\min_{F\in \mathcal{F}}\chi(F)\geq 3$.
One can find that $\mathcal{M}(\mathcal{F})$ always contains some bipartite graphs.
Set
            $$q(\mathcal{F})=\min\{q(F)~|~F\in \mathcal{M}(\mathcal{F})~~\text{is bipartite} \},$$
where $q(F)$ denotes the minimum number of vertices in an independent covering of $F$.

The following lemma is given by Chv\'{a}tal and Hanson \cite{CH-1976}.

\begin{lem} \label{lem3.1}\emph{(\cite{CH-1976})}
Let $n,\nu,\Delta$ be positive integers with $n\geq 2\nu+1$.
\begin{enumerate}[{\rm (i)}]\setlength{\itemsep}{0pt}
\item If $\Delta\leq 2\nu$ and $n\leq 2\nu+\big\lfloor\frac{\nu}{\lfloor(\Delta+1)/2\rfloor}\big\rfloor$, then
$$f(\nu,\Delta)=
   \begin{cases}
    \min\big\{\big\lfloor\frac{n\Delta}{2}\big\rfloor,
    \nu\Delta+\frac{\Delta-1}{2}\big\lfloor\frac{2(n-\nu)}{\Delta+3}\big\rfloor\big\},  & \hbox{if $\Delta$ is odd;}\\
    \frac{n\Delta}{2},  & \hbox{if $\Delta$ is even.}
    \end{cases}
$$
\item If $\Delta\leq 2\nu$ and $n\geq 2\nu+\big\lfloor\frac{\nu}{\lfloor(\Delta+1)/2\rfloor}\big\rfloor+1$, then
  $$f(\nu,\Delta)=\nu\Delta+\Big\lfloor\frac{\nu}{\lfloor(\Delta+1)/2\rfloor}\Big\rfloor \Big\lfloor\frac{\Delta}{2}\Big\rfloor.$$
\item If $\Delta\geq 2\nu+1$, then
$$f(\nu,\Delta)=
   \begin{cases}
    \max\big\{\binom{2\nu+1}{2},\big\lfloor\frac{\nu(n+\Delta-\nu)}{2}\big\rfloor\big\},
    & \hbox{if $n\leq \nu+\Delta$;}\\
    \nu\Delta,  & \hbox{if $n\geq \nu+\Delta+1$.}
    \end{cases}
$$
\end{enumerate}
\end{lem}

For convenience, we define $H(n,k,i)=(M_{k-1}\cup K_1)\nabla T_{n-2k+1,i}$  for $i\in [2]$, where $T_{n-2k+1,1}=E_{n-2k+1}$,
and set $h(n,k,i)=e(H(n,k,i))$.
We now turn our attention to the Tur\'{a}n number of $\mathcal{M}(F^{o})$,
which will enable us to determine that of $F^{o}$.

\begin{lem}\label{lem2.2}
Let $n$ be sufficiently large,
and $F=K_1\nabla F^{\bullet}$ where each component of $F^{\bullet}$ is either a non-trivial tree or an even cycle. Then
 $${\rm EX}(n,\mathcal{M}(F^{o}))=
   \begin{cases}
   \{H(n,\frac{1}{2}{|F^{\bullet}|},1)\},
       & \hbox{if each component of $F^{\bullet}$ is an even cycle;} \\
     \{K_{e(F^{\bullet})}\nabla E_{n-e(F^{\bullet})}\},
       & \hbox{if each component of $F^{\bullet}$ is an edge;} \\
          \{K_{e(F^{\bullet}),n-e(F^{\bullet})}\},
            & \hbox{otherwise.}
              \end{cases}
$$
\end{lem}

\begin{proof}
In the graph $F$, let $w=V(K_1)$,
and let $A=\{v_1,\dots,v_{|A|}\}$ and $B=\{u_1,\dots,u_{|B|}\}$
be two color classes of $F^{\bullet}$ such that $|A|$ is as small as possible.
Let $G$ be an extremal graph for $\mathcal{M}(F^{o})$.
For convenience, we write $q=q(\mathcal{M}(F^{o}))$.
From the definition of $q(\mathcal{M}(F^{o}))$
we know that $K_{q-1,n-q+1}$ is $\mathcal{M}(F^{o})$-free.
Thus,
\begin{align}\label{equ-2A}
e(G)\geq e(K_{q-1,n-q+1})=(q-1)n-(q-1)^2.
\end{align}

Now we present three claims.

\begin{claim}\label{claim-2.0}
$q=e(F^{\bullet})$ if each component of $F^{\bullet}$ is an even cycle,
and  $q=e(F^{\bullet})+1$ otherwise.
\end{claim}
\renewcommand\proofname{\bf Proof of Claim \ref{claim-2.0}}
\begin{proof}
Let $J\in \mathcal{M}(F^{o})$ be a bipartite graph with $q=q(J)$.
By Lemma \ref{lem2.1}, there is an independent set $U$ of $F$ such that $J\in \mathcal{C}(F,U)$.
Suppose first that $w\in U$.
Since $U$ is an independent set of $F$, it holds that $U=\{w\}$.
Let $J(w)$ be the subgraph of $J$ consisting of all edges in $\{e(w,u)~|~u\in V(F^{\bullet})\}$.
Clearly, $J(w)\cong M_{|F^{\bullet}|}$.
It follows that $q(J)\geq q(J(w))=|F^{\bullet}|$.
Suppose then that $w\notin U$.
Apparently, $U\subseteq V(F^{\bullet})$.
Clearly, $U$ is a covering set of $F^{\bullet}$.
Otherwise, $F^{\bullet}-U$ contains an edge $xy$.
Then $wxy$ is a triangle of $J$, which contradicts that $J$ is bipartite.
Let $J^*\in \mathcal{C}(F,U)$, where each edge of $F[U,V(F)\setminus U]$ is of Type I.
Since $U$ is an independent set and a covering set of $F^{\bullet}$,
it follows that
$J^*$ is a bipartite graph with the two partite sets
containing $e(F^{\bullet})+1$ and $|F^{\bullet}|$ vertices respectively.
It is not hard to observe that $q(J)\geq q(J^*)=\min\{e(F^{\bullet})+1,|F^{\bullet}|\}$.
Based on the above discussion,
by the definitions of $q$ and $J$, we get
\begin{align}\label{equ-1A}
q=\min\{q(F)~|~F\in \mathcal{M}(\mathcal{F})~~\text{is bipartite} \}\geq \min\{e(F^{\bullet})+1,|F^{\bullet}|\}.
\end{align}

Now, we define $J_1,J_2,J_3$ and $J_4$ as follows
(see $J_1,J_2,J_3$ and $J_4$ for $F=W_{2k+1}$ in Figure \ref{fig-2.2}):
\begin{itemize}\setlength{\itemsep}{0pt}
\item $J_1\in \mathcal{C}(F,B)$, where each edge of $F[\{w\}\cup A,B]$ is of Type I;
\item $J_2\in \mathcal{C}(F,B)$, where each edge of $F[\{w\}\cup A,B]$ is of Type II.
Clearly, $J_2\cong S_{|A|+1}\cup M_{|B|+e(F^{\bullet})}$;
\item $J_3\in \mathcal{C}(F,B)$, where each edge of $F[\{w\},B]$ is of Type I and each edge of $F[A,B]$ is of Type II.
Clearly, $J_3\cong S_{|F^{\bullet}|+1}\cup M_{e(F^{\bullet})}$;
\item $J_4\in \mathcal{C}(F,\{w\})$, where each edge of $F[A,\{w\}]$ is of Type I and each edge of $F[B,\{w\}]$ is of Type II. Clearly, $q(J_4)=|F^{\bullet}|$, and $J_4$ contains exactly $|B|$ edge components.
\end{itemize}

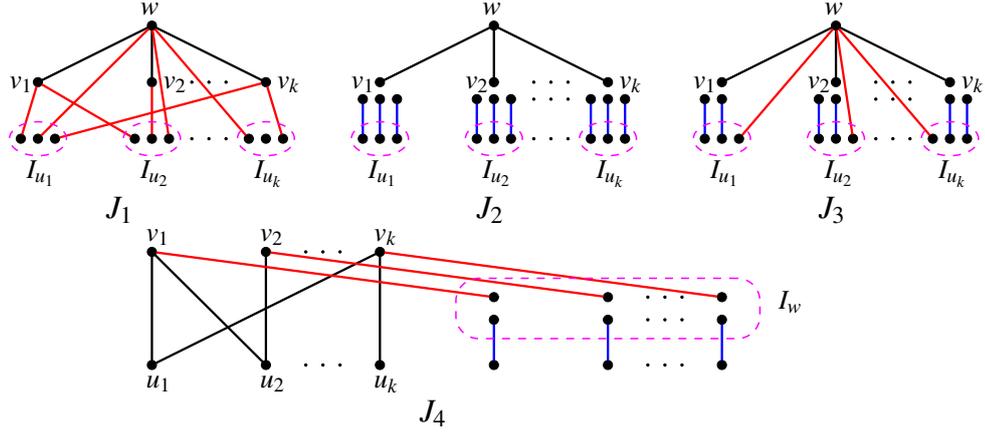
\begin{figure}
\centering
\begin{tikzpicture}[scale=0.75, x=1.00mm, y=1.00mm, inner xsep=0pt, inner ysep=0pt, outer xsep=0pt, outer ysep=0pt]
\definecolor{F}{rgb}{0,0,0}

\node[circle,fill=F,draw,inner sep=0pt,minimum size=1.2mm] (w) at (10,30) {};
\draw(8,32) node[anchor=base west]{\fontsize{10.23}{17.07}\selectfont $w$};

\node[circle,fill=F,draw,inner sep=0pt,minimum size=0.3mm] () at (17,20) {};
\node[circle,fill=F,draw,inner sep=0pt,minimum size=0.3mm] () at (20,20) {};
\node[circle,fill=F,draw,inner sep=0pt,minimum size=0.3mm] () at (23,20) {};

\node[circle,fill=F,draw,inner sep=0pt,minimum size=1.2mm] (u1) at (-10,20) {};
\draw(-15,19) node[anchor=base west]{\fontsize{10.23}{17.07}\selectfont $v_1$};
\node[circle,fill=F,draw,inner sep=0pt,minimum size=1.2mm] (u2) at (10,20) {};
\draw(12,19) node[anchor=base west]{\fontsize{10.23}{17.07}\selectfont $v_2$};

\node[circle,fill=F,draw,inner sep=0pt,minimum size=1.2mm] (uk) at (30,20) {};
\draw(32,19) node[anchor=base west]{\fontsize{10.23}{17.07}\selectfont $v_{k}$};

\node[circle,fill=F,draw,inner sep=0pt,minimum size=1.2mm] (v1a) at (-13,10) {};
\node[circle,fill=F,draw,inner sep=0pt,minimum size=1.2mm] (v1b) at (-10,10) {};
\node[circle,fill=F,draw,inner sep=0pt,minimum size=1.2mm] (v1) at (-7,10) {};

\node[circle,fill=F,draw,inner sep=0pt,minimum size=1.2mm] (v2a) at (7,10) {};
\node[circle,fill=F,draw,inner sep=0pt,minimum size=1.2mm] (v2b) at (10,10) {};
\node[circle,fill=F,draw,inner sep=0pt,minimum size=1.2mm] (v2) at (13,10) {};

\node[circle,fill=F,draw,inner sep=0pt,minimum size=0.3mm] () at (17,10) {};
\node[circle,fill=F,draw,inner sep=0pt,minimum size=0.3mm] () at (20,10) {};
\node[circle,fill=F,draw,inner sep=0pt,minimum size=0.3mm] () at (23,10) {};

\node[circle,fill=F,draw,inner sep=0pt,minimum size=1.2mm] (vka) at (33,10) {};
\node[circle,fill=F,draw,inner sep=0pt,minimum size=1.2mm] (vkb) at (30,10) {};
\node[circle,fill=F,draw,inner sep=0pt,minimum size=1.2mm] (vk) at (27,10) {};

\definecolor{L}{rgb}{0,0,0}
\path[line width=0.3mm, draw=L] (w) -- (u1);
\path[line width=0.3mm, draw=L] (w) -- (u2);
\path[line width=0.3mm, draw=L] (w) -- (uk);

\definecolor{L}{rgb}{1,0,0}

\path[line width=0.3mm, draw=L] (w) -- (v1b);

\path[line width=0.3mm, draw=L] (u1) -- (v1a);
\path[line width=0.3mm, draw=L] (uk) -- (v1);
\path[line width=0.3mm, draw=L] (u1) -- (v2a);
\path[line width=0.3mm, draw=L] (u2) -- (v2b);
\path[line width=0.3mm, draw=L] (w) -- (v2);

\path[line width=0.3mm, draw=L] (uk) -- (vka);
\path[line width=0.3mm, draw=L] (w) -- (vk);


\node[circle,fill=F,draw,inner sep=0pt,minimum size=1.2mm] (w) at (70,30) {};
\draw(68,32) node[anchor=base west]{\fontsize{10.23}{17.07}\selectfont $w$};

\node[circle,fill=F,draw,inner sep=0pt,minimum size=0.3mm] () at (77,20) {};
\node[circle,fill=F,draw,inner sep=0pt,minimum size=0.3mm] () at (80,20) {};
\node[circle,fill=F,draw,inner sep=0pt,minimum size=0.3mm] () at (83,20) {};

\node[circle,fill=F,draw,inner sep=0pt,minimum size=1.2mm] (u1) at (50,20) {};
\draw(45,19) node[anchor=base west]{\fontsize{10.23}{17.07}\selectfont $v_1$};
\node[circle,fill=F,draw,inner sep=0pt,minimum size=1.2mm] (u2) at (70,20) {};
\draw(65,19) node[anchor=base west]{\fontsize{10.23}{17.07}\selectfont $v_2$};

\node[circle,fill=F,draw,inner sep=0pt,minimum size=1.2mm] (uk) at (90,20) {};
\draw(92,19) node[anchor=base west]{\fontsize{10.23}{17.07}\selectfont $v_{k}$};

\node[circle,fill=F,draw,inner sep=0pt,minimum size=1.2mm] (v1a) at (47,10) {};
\node[circle,fill=F,draw,inner sep=0pt,minimum size=1.2mm] (v1b) at (50,10) {};
\node[circle,fill=F,draw,inner sep=0pt,minimum size=1.2mm] (v1) at (53,10) {};

\node[circle,fill=F,draw,inner sep=0pt,minimum size=1.2mm] (v2a) at (67,10) {};
\node[circle,fill=F,draw,inner sep=0pt,minimum size=1.2mm] (v2b) at (70,10) {};
\node[circle,fill=F,draw,inner sep=0pt,minimum size=1.2mm] (v2) at (73,10) {};

\node[circle,fill=F,draw,inner sep=0pt,minimum size=0.3mm] () at (77,10) {};
\node[circle,fill=F,draw,inner sep=0pt,minimum size=0.3mm] () at (80,10) {};
\node[circle,fill=F,draw,inner sep=0pt,minimum size=0.3mm] () at (83,10) {};

\node[circle,fill=F,draw,inner sep=0pt,minimum size=1.2mm] (vka) at (93,10) {};
\node[circle,fill=F,draw,inner sep=0pt,minimum size=1.2mm] (vkb) at (90,10) {};
\node[circle,fill=F,draw,inner sep=0pt,minimum size=1.2mm] (vk) at (87,10) {};

\node[circle,fill=F,draw,inner sep=0pt,minimum size=1.2mm] (w1a) at (47,17) {};
\node[circle,fill=F,draw,inner sep=0pt,minimum size=1.2mm] (w1b) at (50,17) {};
\node[circle,fill=F,draw,inner sep=0pt,minimum size=1.2mm] (w1) at (53,17) {};

\node[circle,fill=F,draw,inner sep=0pt,minimum size=1.2mm] (w2a) at (67,17) {};
\node[circle,fill=F,draw,inner sep=0pt,minimum size=1.2mm] (w2b) at (70,17) {};
\node[circle,fill=F,draw,inner sep=0pt,minimum size=1.2mm] (w2) at (73,17) {};

\node[circle,fill=F,draw,inner sep=0pt,minimum size=0.3mm] () at (77,17) {};
\node[circle,fill=F,draw,inner sep=0pt,minimum size=0.3mm] () at (80,17) {};
\node[circle,fill=F,draw,inner sep=0pt,minimum size=0.3mm] () at (83,17) {};

\node[circle,fill=F,draw,inner sep=0pt,minimum size=1.2mm] (wka) at (93,17) {};
\node[circle,fill=F,draw,inner sep=0pt,minimum size=1.2mm] (wkb) at (90,17) {};
\node[circle,fill=F,draw,inner sep=0pt,minimum size=1.2mm] (wk) at (87,17) {};

\definecolor{L}{rgb}{0,0,0}
\path[line width=0.3mm, draw=L] (w) -- (u1);
\path[line width=0.3mm, draw=L] (w) -- (u2);
\path[line width=0.3mm, draw=L] (w) -- (uk);

\definecolor{L}{rgb}{0,0,1}

\path[line width=0.3mm, draw=L] (w1) -- (v1);
\path[line width=0.3mm, draw=L] (w2) -- (v2);
\path[line width=0.3mm, draw=L] (wk) -- (vk);

\path[line width=0.3mm, draw=L] (w1a) -- (v1a);
\path[line width=0.3mm, draw=L] (w1b) -- (v1b);
\path[line width=0.3mm, draw=L] (w2a) -- (v2a);
\path[line width=0.3mm, draw=L] (w2b) -- (v2b);
\path[line width=0.3mm, draw=L] (wka) -- (vka);
\path[line width=0.3mm, draw=L] (wkb) -- (vkb);

\node[circle,fill=F,draw,inner sep=0pt,minimum size=1.2mm] (w) at (130,30) {};
\draw(128,32) node[anchor=base west]{\fontsize{10.23}{17.07}\selectfont $w$};

\node[circle,fill=F,draw,inner sep=0pt,minimum size=0.3mm] () at (137,20) {};
\node[circle,fill=F,draw,inner sep=0pt,minimum size=0.3mm] () at (140,20) {};
\node[circle,fill=F,draw,inner sep=0pt,minimum size=0.3mm] () at (143,20) {};

\node[circle,fill=F,draw,inner sep=0pt,minimum size=1.2mm] (u1) at (110,20) {};
\draw(105,19) node[anchor=base west]{\fontsize{10.23}{17.07}\selectfont $v_1$};
\node[circle,fill=F,draw,inner sep=0pt,minimum size=1.2mm] (u2) at (130,20) {};
\draw(125,19) node[anchor=base west]{\fontsize{10.23}{17.07}\selectfont $v_2$};

\node[circle,fill=F,draw,inner sep=0pt,minimum size=1.2mm] (uk) at (150,20) {};
\draw(152,19) node[anchor=base west]{\fontsize{10.23}{17.07}\selectfont $v_{k}$};

\node[circle,fill=F,draw,inner sep=0pt,minimum size=1.2mm] (v1a) at (107,10) {};
\node[circle,fill=F,draw,inner sep=0pt,minimum size=1.2mm] (v1b) at (110,10) {};
\node[circle,fill=F,draw,inner sep=0pt,minimum size=1.2mm] (v1) at (113,10) {};

\node[circle,fill=F,draw,inner sep=0pt,minimum size=1.2mm] (v2a) at (127,10) {};
\node[circle,fill=F,draw,inner sep=0pt,minimum size=1.2mm] (v2b) at (130,10) {};
\node[circle,fill=F,draw,inner sep=0pt,minimum size=1.2mm] (v2) at (133,10) {};

\node[circle,fill=F,draw,inner sep=0pt,minimum size=0.3mm] () at (137,10) {};
\node[circle,fill=F,draw,inner sep=0pt,minimum size=0.3mm] () at (140,10) {};
\node[circle,fill=F,draw,inner sep=0pt,minimum size=0.3mm] () at (143,10) {};

\node[circle,fill=F,draw,inner sep=0pt,minimum size=1.2mm] (vka) at (153,10) {};
\node[circle,fill=F,draw,inner sep=0pt,minimum size=1.2mm] (vkb) at (150,10) {};
\node[circle,fill=F,draw,inner sep=0pt,minimum size=1.2mm] (vk) at (147,10) {};

\node[circle,fill=F,draw,inner sep=0pt,minimum size=1.2mm] (w1a) at (107,17) {};
\node[circle,fill=F,draw,inner sep=0pt,minimum size=1.2mm] (w1b) at (110,17) {};

\node[circle,fill=F,draw,inner sep=0pt,minimum size=1.2mm] (w2a) at (127,17) {};
\node[circle,fill=F,draw,inner sep=0pt,minimum size=1.2mm] (w2b) at (130,17) {};

\node[circle,fill=F,draw,inner sep=0pt,minimum size=0.3mm] () at (137,17) {};
\node[circle,fill=F,draw,inner sep=0pt,minimum size=0.3mm] () at (140,17) {};
\node[circle,fill=F,draw,inner sep=0pt,minimum size=0.3mm] () at (143,17) {};

\node[circle,fill=F,draw,inner sep=0pt,minimum size=1.2mm] (wka) at (153,17) {};
\node[circle,fill=F,draw,inner sep=0pt,minimum size=1.2mm] (wkb) at (150,17) {};

\definecolor{L}{rgb}{0,0,0}
\path[line width=0.3mm, draw=L] (w) -- (u1);
\path[line width=0.3mm, draw=L] (w) -- (u2);
\path[line width=0.3mm, draw=L] (w) -- (uk);

\definecolor{L}{rgb}{1,0,0}

\path[line width=0.3mm, draw=L] (w) -- (v1);
\path[line width=0.3mm, draw=L] (w) -- (v2);
\path[line width=0.3mm, draw=L] (w) -- (vk);

\definecolor{L}{rgb}{0,0,1}

\path[line width=0.3mm, draw=L] (w1a) -- (v1a);
\path[line width=0.3mm, draw=L] (w1b) -- (v1b);
\path[line width=0.3mm, draw=L] (w2a) -- (v2a);
\path[line width=0.3mm, draw=L] (w2b) -- (v2b);
\path[line width=0.3mm, draw=L] (wka) -- (vka);
\path[line width=0.3mm, draw=L] (wkb) -- (vkb);

\definecolor{L}{rgb}{1,0,1}

\path[line width=0.20mm, draw=L,dashed] (-10,10) ellipse (5.00mm and 3.00mm);
\path[line width=0.20mm, draw=L,dashed] (10,10) ellipse (5.00mm and 3.00mm);
\path[line width=0.20mm, draw=L,dashed] (30,10) ellipse (5.00mm and 3.00mm);
\path[line width=0.20mm, draw=L,dashed] (50,10) ellipse (5.00mm and 3.00mm);
\path[line width=0.20mm, draw=L,dashed] (70,10) ellipse (5.00mm and 3.00mm);
\path[line width=0.20mm, draw=L,dashed] (90,10) ellipse (5.00mm and 3.00mm);
\path[line width=0.20mm, draw=L,dashed] (110,10) ellipse (5.00mm and 3.00mm);
\path[line width=0.20mm, draw=L,dashed] (130,10) ellipse (5.00mm and 3.00mm);
\path[line width=0.20mm, draw=L,dashed] (150,10) ellipse (5.00mm and 3.00mm);

\draw(2,-4) node[anchor=base west]{\fontsize{12.23}{17.07}\selectfont $J_1$};
\draw(67,-4) node[anchor=base west]{\fontsize{12.23}{17.07}\selectfont $J_2$};
\draw(127,-4) node[anchor=base west]{\fontsize{12.23}{17.07}\selectfont $J_3$};

\draw(-12,3) node[anchor=base west]{\fontsize{10.23}{17.07}\selectfont $I_{u_1}$};
\draw(8,3) node[anchor=base west]{\fontsize{10.23}{17.07}\selectfont $I_{u_2}$};
\draw(28,3) node[anchor=base west]{\fontsize{10.23}{17.07}\selectfont $I_{u_k}$};

\draw(48,3) node[anchor=base west]{\fontsize{10.23}{17.07}\selectfont $I_{u_1}$};
\draw(68,3) node[anchor=base west]{\fontsize{10.23}{17.07}\selectfont $I_{u_2}$};
\draw(88,3) node[anchor=base west]{\fontsize{10.23}{17.07}\selectfont $I_{u_k}$};

\draw(108,3) node[anchor=base west]{\fontsize{10.23}{17.07}\selectfont $I_{u_1}$};
\draw(128,3) node[anchor=base west]{\fontsize{10.23}{17.07}\selectfont $I_{u_2}$};
\draw(148,3) node[anchor=base west]{\fontsize{10.23}{17.07}\selectfont $I_{u_k}$};


\definecolor{F}{rgb}{0,0,0}
\node[circle,fill=F,draw,inner sep=0pt,minimum size=1.2mm] (u1) at (10,-10) {};
\draw(9,-8) node[anchor=base west]{\fontsize{10.23}{17.07}\selectfont $v_1$};
\node[circle,fill=F,draw,inner sep=0pt,minimum size=1.2mm] (u2) at (30,-10) {};
\draw(29,-8) node[anchor=base west]{\fontsize{10.23}{17.07}\selectfont $v_2$};

\node[circle,fill=F,draw,inner sep=0pt,minimum size=0.3mm] () at (37,-10) {};
\node[circle,fill=F,draw,inner sep=0pt,minimum size=0.3mm] () at (40,-10) {};
\node[circle,fill=F,draw,inner sep=0pt,minimum size=0.3mm] () at (43,-10) {};

\node[circle,fill=F,draw,inner sep=0pt,minimum size=1.2mm] (uk) at (50,-10) {};
\draw(49,-8) node[anchor=base west]{\fontsize{10.23}{17.07}\selectfont $v_{k}$};

\node[circle,fill=F,draw,inner sep=0pt,minimum size=1.2mm] (v1) at (10,-30) {};
\draw(9,-34) node[anchor=base west]{\fontsize{10.23}{17.07}\selectfont $u_1$};
\node[circle,fill=F,draw,inner sep=0pt,minimum size=1.2mm] (v2) at (30,-30) {};
\draw(29,-34) node[anchor=base west]{\fontsize{10.23}{17.07}\selectfont $u_2$};
\node[circle,fill=F,draw,inner sep=0pt,minimum size=0.3mm] () at (37,-30) {};
\node[circle,fill=F,draw,inner sep=0pt,minimum size=0.3mm] () at (40,-30) {};
\node[circle,fill=F,draw,inner sep=0pt,minimum size=0.3mm] () at (43,-30) {};

\node[circle,fill=F,draw,inner sep=0pt,minimum size=1.2mm] (vk) at (50,-30) {};
\draw(49,-34) node[anchor=base west]{\fontsize{10.23}{17.07}\selectfont $u_k$};

\definecolor{L}{rgb}{0,0,0}
\path[line width=0.3mm, draw=L] (v1) -- (u1);
\path[line width=0.3mm, draw=L] (v2) -- (u2);
\path[line width=0.3mm, draw=L] (vk) -- (uk);
\path[line width=0.3mm, draw=L] (v2) -- (u1);
\path[line width=0.3mm, draw=L] (v1) -- (uk);

\node[circle,fill=F,draw,inner sep=0pt,minimum size=1.2mm] (w1) at (70,-18) {};
\node[circle,fill=F,draw,inner sep=0pt,minimum size=1.2mm] (w2) at (90,-18) {};
\node[circle,fill=F,draw,inner sep=0pt,minimum size=0.3mm] () at (97,-18) {};
\node[circle,fill=F,draw,inner sep=0pt,minimum size=0.3mm] () at (100,-18) {};
\node[circle,fill=F,draw,inner sep=0pt,minimum size=0.3mm] () at (103,-18) {};

\node[circle,fill=F,draw,inner sep=0pt,minimum size=1.2mm] (wk) at (110,-18) {};

\definecolor{L}{rgb}{1,0,0}
\path[line width=0.3mm, draw=L] (u1) -- (w1);
\path[line width=0.3mm, draw=L] (u2) -- (w2);
\path[line width=0.3mm, draw=L] (uk) -- (wk);

\node[circle,fill=F,draw,inner sep=0pt,minimum size=1.2mm] (w1L) at (70,-22) {};
\node[circle,fill=F,draw,inner sep=0pt,minimum size=1.2mm] (w2L) at (90,-22) {};
\node[circle,fill=F,draw,inner sep=0pt,minimum size=0.3mm] () at (97,-22) {};
\node[circle,fill=F,draw,inner sep=0pt,minimum size=0.3mm] () at (100,-22) {};
\node[circle,fill=F,draw,inner sep=0pt,minimum size=0.3mm] () at (103,-22) {};

\node[circle,fill=F,draw,inner sep=0pt,minimum size=1.2mm] (wkL) at (110,-22) {};

\node[circle,fill=F,draw,inner sep=0pt,minimum size=1.2mm] (w1R) at (70,-30) {};
\node[circle,fill=F,draw,inner sep=0pt,minimum size=1.2mm] (w2R) at (90,-30) {};
\node[circle,fill=F,draw,inner sep=0pt,minimum size=0.3mm] () at (97,-30) {};
\node[circle,fill=F,draw,inner sep=0pt,minimum size=0.3mm] () at (100,-30) {};
\node[circle,fill=F,draw,inner sep=0pt,minimum size=0.3mm] () at (103,-30) {};

\node[circle,fill=F,draw,inner sep=0pt,minimum size=1.2mm] (wkR) at (110,-30) {};

\definecolor{L}{rgb}{0,0,1}
\path[line width=0.3mm, draw=L] (w1L) -- (w1R);
\path[line width=0.3mm, draw=L] (w2L) -- (w2R);
\path[line width=0.3mm, draw=L] (wkL) -- (wkR);

\draw(120,-20) node[anchor=base west]{\fontsize{10.23}{17.07}\selectfont $I_w$};
\draw(57,-40) node[anchor=base west]{\fontsize{12.23}{17.07}\selectfont $J_4$};

\definecolor{L}{rgb}{1,0,1}
\node[draw=L, line width=0.2mm,dashed, minimum width=40mm, minimum height=8mm, rounded corners=3mm] (rect2b) at (90,-20) {};

\end{tikzpicture}
\caption{The graphs $J_1,J_2,J_3,J_4\in \mathcal{C}(W_{2k+1})$.}{\label{fig-2.2}}
\end{figure}

We first consider the case that each component of $F^{\bullet}$ is an even cycle.
Then, $q(J_4)=|F^{\bullet}|$ and $e(F^{\bullet})=|F^{\bullet}|$.
Combining these with \eqref{equ-1A}, we obtain $q=e(F^{\bullet})$.
We then consider the case that $F^{\bullet}$ contains tree components.
Then, $q(J_3)=e(F^{\bullet})+1$ and $e(F^{\bullet})+1\leq |F^{\bullet}|$.
Combining these with \eqref{equ-1A}, we obtain $q=e(F^{\bullet})+1$.
This completes the proof of Claim \ref{claim-2.0}.
\end{proof}

\begin{claim}\label{claim-2.1}
Let $R\subseteq V(G)$ be the set of vertices of degree at least $80|F^{\bullet}|$.
Then $|R|=q-1$.
\end{claim}

\renewcommand\proofname{\bf Proof of Claim \ref{claim-2.1}}
\begin{proof}
Let $u^*$ be a vertex of $G$ with maximum degree and set $X=\{u^*\}\cup N_G(u^*)$.
Since $n$ is sufficiently large,
the Handshaking Lemma and \eqref{equ-2A} yields $d_G(u^*)\geq \lceil\frac{2e(G)}{n}\rceil=2(q-1)$.
By Claim \ref{claim-2.0}, we get $q\geq e(F^{\bullet})$.
Since $|A|$ is as small as possible, we see that $e(F^{\bullet})\geq |A|$.
These yield that $d_G(u^*)\geq  2(q-1)\geq |A|$.
Consequently, $G[X]$ contains a copy of $S_{|A|+1}$.

We first show that $R\neq \varnothing$. Suppose for contradiction that $R = \varnothing$.
The definition of $R$ then implies that $\Delta(G) < 80|F^{\bullet}|$.
Observe that $J_2\cong S_{|A|+1}\cup M_{|B|+e(F^{\bullet})}$. Since $G[X]$ already contains a copy of $S_{|A|+1}$, the fact that $G$ is $J_2$-free forces the subgraph $G-X$ to be $M_{|B|+e(F^{\bullet})}$-free.
Applying Lemma \ref{lem3.1} to $G-X$, we obtain
$$e(G-X)\leq f(|B|+e(F^{\bullet})-1,80|F^{\bullet}|-1)\leq (|B|+e(F^{\bullet})-1)80|F^{\bullet}|.$$
Since $\Delta(G)< 80|F^{\bullet}|$,
we clearly have $|X|=d_G(u^*)+1\leq 80|F^{\bullet}|$.
Consequently,
\begin{align*}
e(G) \leq e(G-X)+\sum_{v\in X}d_G(v) \leq (|B|+e(F^{\bullet})-1)80|F^{\bullet}|+(80|F^{\bullet}|-1)80|F^{\bullet}|,
\end{align*}
which contradicts \eqref{equ-2A} since $n$ is sufficiently large.
This contradiction establishes that $R \neq \emptyset$. Furthermore, since $u^*$ is a vertex of maximum degree in $G$, we have $u^* \in R$.

We next show that $G-R$ is $M_{e(F^{\bullet})}$-free.
Suppose, to the contrary, that $G-R$ has a subgraph $H'$ isomorphic to $M_{e(F^{\bullet})}$.
Since $u^*\in R$, we get
$$d_G(u^*)\geq 80|F^{\bullet}|\geq |F^{\bullet}|+2e(F^{\bullet})= |F^{\bullet}|+|H'|.$$
This implies that $u^*$ has at least $|F^{\bullet}|$ neighbors outside $V(H')$, and thus $G-V(H')$ contains a star $S_{|F^{\bullet}|+1}$ centered at $u^*$.
It follows that $G$ contains a subgraph $S_{|F^{\bullet}|+1}\cup H'$ isomorphic to $J_3$,
which leads to a contradiction.
Hence, $G-R$ must be $M_{e(F^{\bullet})}$-free.
By the definition of $R$, we have $\Delta(G-R)\leq 80|F^{\bullet}|-1$.
Applying Lemma \ref{lem3.1} yields
\begin{align}\label{equ-3A}
e(G-R)\leq f(e(F^{\bullet})-1,80|F^{\bullet}|-1)\leq (e(F^{\bullet})-1)80|F^{\bullet}|.
\end{align}

We now bound the size of $R$.
First, if $|R|\leq q-2$, then
\begin{align*}
e(G) \leq e(G-R)+\sum_{v\in R}d_G(v) \leq (e(F^{\bullet})-1)80|F^{\bullet}|+(q-2)n,
\end{align*}
which contradicts \eqref{equ-2A}.

Next, suppose for contradiction that $|R|\geq q+1$.
By Claim \ref{claim-2.0}, we get $|R|-1\geq q\geq e(F^{\bullet})$.
Since every vertex in $R$ has degree at least $80|F^{\bullet}|$, the bipartite subgraph $G[R \setminus \{u^*\}, V(G) \setminus R]$ contains a matching $H' \cong M_{e(F^{\bullet})}$.
Since $d_G(u^*)\geq 80|F^{\bullet}|\geq |F^{\bullet}|+|H'|$,
the vertex $u^*$ must have at least $|F^{\bullet}|$ neighbors in $V(G) \setminus V(H')$.
Thus, $G-V(H')$ contains a star $S_{|F^{\bullet}|+1}$ centered at $u_1$.
Consequently, $G$ contains $H' \cup S_{|F^{\bullet}|+1} \cong J_3$ as a subgraph, leading to a contradiction.
Therefore, we must have $|R| \in \{q-1, q\}$.

It remains to rule out the case $|R| = q$.
Suppose for contradiction that $|R| = q$.
We distinguish two cases depending on the structure of $F^{\bullet}$:

\medskip
\noindent\textbf{Case 1.} $F^{\bullet}$ contains tree components.

By Claim \ref{claim-2.0}, we get $q=e(F^{\bullet})+1$.
Since $|R| = q = e(F^{\bullet}) + 1$, an argument analogous to the one above shows that $G$ contains a copy of $J_3 \cong S_{|F^{\bullet}|+1} \cup M_{e(F^{\bullet})}$, a contradiction.

\medskip
\noindent\textbf{Case 2.} $F^{\bullet}$ is the disjoint union of even cycles.

Specifically, let $F^{\bullet}\cong \cup_{i=1}^{j}C_{2k_i}$ and set $k=\sum_{i=1}^{j}k_i$,
which implies that $q=2k$.
Let the vertices of $R$ be ordered as $R = \{w_1, w_2, \dots, w_{2k}\}$ such that $d_G(w_1) \ge d_G(w_2) \ge \dots \ge d_G(w_{2k})$.
If $d_G(w_{k-1})+d_G(w_k)\leq 1.1n$,
then $d_G(w_i)\leq d_G(w_k)\leq 0.55n$ for all $i\in \{k+1,\dots,2k\}$.
Combining these with \eqref{equ-3A}, we obtain that
\begin{align*}
e(G) &\leq e(G-R)+\sum_{w\in R}d_G(w)
\leq (2k-1)160k+(k-2)n+1.1n+k\times 0.55n\\
&<(1.55k-0.8)n.
\end{align*}
Since $k\geq 2$, we have $1.55k-0.8<2k-1.1$,
yielding $e(G)<(2k-1.1)n$,
which contradicts \eqref{equ-2A}.
Hence, we must have $d_G(w_{k-1})+d_G(w_k)>1.1n$.
For any distinct $i_1, i_2 \in [k]$, this inequality implies that
 $$|N_G(w_{i_1})\cap N_G(w_{i_2})|\geq |N_G(w_{i_1})|+|N_G(w_{i_2})|-n\geq d_G(w_{k-1})+d_G(w_k)-n\geq 0.1n.$$
By greedily selecting common neighbors, we find that $G$ contains vertex-disjoint cycles $C^1, C^2, \dots, C^j$ such that each $C^i$ contains exactly $k_i$ vertices from $\{w_1, \dots, w_k\}$.
Given that $|R|=q=2k$ and each vertex in $R$ has degree at least $80|F^{\bullet}|$,
it follows that $G$ contains a copy of $J_4$, a contradiction.

By the above case analysis, we conclude that $|R|=q-1$, as desired.
\end{proof}

\begin{claim}\label{claim-2.2}
(a) $G-R$ is empty.
(b) $G[R,V(G)\setminus R]$ is a complete bipartite graph.
\end{claim}

\renewcommand\proofname{\bf Proof of Claim \ref{claim-2.2}}
\begin{proof}
In view of Claim \ref{claim-2.1}, we have $|R|=q-1$.
Let us label the vertices of $R$ as $\{w_1, w_2, \dots, w_{q-1}\}$ such that their degrees are in non-increasing order, i.e., $d_{G}(w_1) \geq d_{G}(w_2) \geq \cdots \geq d_{G}(w_{q-1})$.
Applying \eqref{equ-3A}, we can bound the number of edges in $G$ by
\begin{align*}
e(G) \leq e(G-R)+\sum_{v\in R}d_G(v) \leq (e(F^{\bullet})-1)80|F^{\bullet}|+(q-2)n+d_{G}(w_{q-1}).
\end{align*}
Combining this bound with \eqref{equ-2A} yields $d_{G}(w_{q-1})\geq \frac{2q-3}{2q-2}n$.
Thus,
\begin{align}\label{equ-4A}
\Big|\bigcap_{w\in R}N_G(w)\Big|\geq |R|\cdot \frac{2q-3}{2q-2}n-(q-2)n=\frac{n}{2}.
\end{align}

($a$) We first consider the case where $F^{\bullet}$ is the disjoint union of even cycles,
that is, $F^{\bullet}\cong \cup_{i=1}^{j}C_{2k_i}$.
Set $k=\sum_{i=1}^{j}k_i$.
By Claims \ref{claim-2.0} and \ref{claim-2.1}, $|R|=2k-1$.
Let $J_4'$ be the graph obtained from $J_4$ by removing one of its edge components.
Then, $J_4\subseteq K_{k,2k}\cup M_k$ and $J_4'\subseteq K_{k,2k}\cup M_{k-1}\subseteq K_{2k-1,3k-1}$.
Suppose that $G-R$ contains at least one edge, say $xy$.
By \eqref{equ-4A}, there exist at least $3k-1$ vertices in $V(G)\setminus (R\cup \{x,y\})$ that are adjacent to all vertices in $R$.
Since $|R|=2k-1$, this common neighborhood structure guarantees that $G-\{x,y\}$ contains a subgraph $J''$ isomorphic to $J_4'$.
Since the edge $xy$ is vertex-disjoint from $J''$, the union $J''\cup \{xy\}$ forms a subgraph of $G$ isomorphic to $J_4$, which contradicts the assumption that $G$ is $\mathcal{M}(F^o)$-free.
Thus, $G-R$ is not empty.

Now we turn to the case where $F^{\bullet}$ contains tree components.
By Claims \ref{claim-2.0} and \ref{claim-2.1}, we have $|R|=e(F^{\bullet})$.
Recall that $J_3\cong S_{|F^{\bullet}|+1}\cup M_{e(F^{\bullet})}$.
Let $J_3'$ be the graph obtained from $J_3$ by removing one of its edge components,
so that $J_3'\cong S_{|F^{\bullet}|+1}\cup M_{e(F^{\bullet})-1}$.
Suppose that $G-R$ contains at least one edge, say $xy$.
By \eqref{equ-4A}, there are at least $|F^{\bullet}|+e(F^{\bullet})-1$ vertices in $V(G)\setminus (R\cup \{x,y\})$ that are adjacent to all vertices in $R$.
Since $|R|=e(F^{\bullet})$, it follows that $G-\{x,y\}$ contains a subgraph $J' \cong J_3'$.
Consequently, $G$ contains a subgraph $J' \cup \{xy\}$ isomorphic to $J_3$, yielding a contradiction.

Therefore, in both cases, $G-R$ is empty.

($b$) It suffices to prove that every vertex in $V(G)\setminus R$ is adjacent to all vertices in $R$.
Suppose to the contrary, then there exists a vertex $v_0\in V(G)\setminus R$ that is non-adjacent to some vertex $u_{i_0}$ in $R$.
Let $G'=G+\{v_0u_{i_0}\}$.
Clearly, $e(G')>e(G)$.
By the choice of $G$, we know that $G'$ contains a member $J\in \mathcal{M}(F^{o})$ and $v_0\in V(J)$.
By \eqref{equ-4A}, there exists a vertex $v'\in V(G)\setminus (R\cup V(J))$ such that $v'$ is adjacent to all vertices in $R$.
Note that $N_{J}(v_0)\subseteq R=N_G(v')$.
Then $G[(V(J)\setminus \{v_0\})\cup \{v'\}]$ contains a copy of $J$, a contradiction.
Therefore, $G[R,V(G)\setminus R]$ is a complete bipartite graph.
\end{proof}

Now, we can divide the proof into the following three cases.

\medskip
\noindent{\textbf{Case 1.} $F^{\bullet}$ is the disjoint union of even cycles.}

Suppose that $F^{\bullet} \cong \bigcup_{i=1}^{j} C_{2k_i}$ and let $k = \sum_{i=1}^{j} k_i$.
By Claims \ref{claim-2.0} and \ref{claim-2.1}, we have $|R| = 2k-1$.
Observe that $J_1$ is a tree with $d_{J_1}(v_i)=3$ for each $i$.
If $G[R]$ contains a copy of $P_3$, then $J_1$ can be embedded into $G$ by setting $E(G[R]) \cap E(J_1) = E(P_3) = \{wv_1, wv_2\}$.
Thus, $G[R]$ must be $P_3$-free, which implies that $G[R] \subseteq M_{k-1} \cup K_1$.
For convenience,
we denote $H=(M_{k-1}\cup K_1)\nabla E_{n-2k+1}$ and set $S=V(M_{k-1}\cup K_1)$.
Combining Claims \ref{claim-2.1} and \ref{claim-2.2},
we deduce that $G$ is a subgraph of $H$.

Our goal is to prove that $G\cong H$.
Then it suffices to show that $H$ is $\mathcal{C}(F)$-free.
Suppose to the contrary that $H$ contains a member $J\in \mathcal{C}(F)$.
Then there exists an independent set $U$ of $F$ such that $J\in \mathcal{C}(F,U)$.
Since $J$ is a subgraph of $H$,
we have
\begin{align}\label{align-9Y}
\beta(J)\leq \beta(H)=|S|=2k-1.
\end{align}
If $w\in U$, then $J$ contains $|F^{\bullet}|=2k$ independent edges associated with $w$.
This implies that $\beta(J)\geq 2k$, which contradicts \eqref{align-9Y}.
Consequently, we must have $w\notin U$ and $U\subseteq V(F^{\bullet})$.

We first show that $w \in S$.
Otherwise, $w\notin S$.
Since $w \notin U$, it follows that $w$ is a vertex of $J$, and hence of $H$.
Let $J(w)$ be the subgraph of $J$ consisting of the edges in
$$\{wu \mid u \in N_F(w) \setminus U\} \cup \{e(w,u) \mid u \in N_F(w) \cap U\}.$$
Clearly, $J(w)\subseteq J \subseteq H$.
Since $w \notin S$ and $V(H) \setminus S$ is an independent set,
it follows that $u \in S$ for every vertex $u \in N_F(w) \setminus U$.
Similarly, at least one endpoint of $e(w,u)$ must belong to $S$ for every vertex $u \in N_F(w) \cap U$.
Furthermore, since $e(J(w)) = 2k$ and these edges (or vertices) are mutually disjoint (except possibly sharing $w$, which is not in $S$), $J(w)$ must contain at least $2k$ distinct vertices from $S$.
This contradiction to $|S| = 2k-1$ completes the proof that $w \in S$.

Since $F^{\bullet}$ is the disjoint union of even cycles and $U$ is an independent set of $F$,
it is  clear that $|U|\leq k$.
Secondly, we prove that $|U|=k$.
Otherwise, $|U|\leq k-1$.
For convenience, we denote $U'=V(F^{\bullet})\setminus U$.
Then $|U'|\geq k+1$.
From the construction of $F$ we know that $F[U']$ is not an empty graph.
Since $J[U']=F[U']$, the induced subgraph $J[U']$ is also not an empty graph.
Because $U'\subseteq N_F(w)$ and $V(F)\cap V(J)=U'\cup \{w\}$, it follows that $U'\subseteq N_J(w)$.
We next claim that $|S \cap U'| \le 1$.
Recall that $J \subseteq H$ and $H[S] = M_{k-1} \cup K_1$.
Since $w \in S$ (as proved above) and $U' \subseteq N_J(w) \subseteq N_H(w)$, any vertex $v \in S \cap U'$ must be adjacent to $w$ in $H[S]$.
However, since the maximum degree of $H[S]$ is at most $1$, the vertex $w$ can have at most one neighbor in $H[S]$.
Thus, we indeed have $|S \cap U'| \le 1$.
If $|S \cap U'| = 0$, then $U'$ is entirely contained in $V(H) \setminus S$.
Since $V(H) \setminus S$ is an independent set,
both $H[U']$ and $F[U']$ must be empty graphs,
which contradicts that $F[U']$ is not an empty graph.
This contradiction establishes that $|S \cap U'| = 1$.

Without loss of generality, we may assume that $S \cap U' = \{u_k\}$.
Since $U'\setminus \{u_k\}\subseteq V(H)\setminus S=V(E_{n-2k+1})$,
the set $U'\setminus \{u_k\}$ must be independent.
Furthermore, the inequality $|U'| \ge k+1$ forces
  $|U'| = k+1$ and $U' \setminus \{u_k\} = \{v_1, \dots, v_k\}$,
and hence $U = \{u_1, \dots, u_{k-1}\}$.
For each $uv\in E(F^{\bullet}[\{u_1,\dots,u_{k-1}\},\{v_1,\dots,v_{k}\}])$,
exactly one endpoint of the edge $e(u,v)$ belongs to $S$.
Since $e(F^{\bullet}[\{u_1,\dots,u_{k-1}\},\{v_1,\dots,v_{k}\}])=2(k-1)$ and $w,u_k\in S$,
it follows that $|S|\geq 2+2(k-1)=2k$, which contradicts \eqref{align-9Y}.
Therefore, we must have $|U|=k$.

Since $|U|=k$ and $U$ is an independent set,
we may assume without loss of generality that $U=\{u_1,\dots,u_k\}$.
As $J \subseteq H$, $H[S]=M_{k-1}\cup K_1$ and $w\in S$,
we obtain that $|N_{J}(w)\cap S|\leq 1$.
This, together with $U'=\{v_1,\dots,v_k\}\subseteq N_{J}(w)$, gives $|U'\cap S|\leq 1$.
If $|U'\cap S|=1$,
we may assume without loss of generality that $v_k\in S$.
Then for each $uv_i\in E(F^{\bullet}[U,\{v_1,\dots,v_{k-1}\}])$,
exactly one endpoint of $e(u,v_i)$ belongs to $S$.
Note that $w,v_k\in S$ and $e(F^{\bullet}[U,\{v_1,\dots,v_{k-1}\}])=2(k-1)$.
Then $|S|\geq 2k$, which contradicts $|S|=2k-1$.
If $|U'\cap S|=0$,
then $U'\subseteq V(H)\setminus S$.
Moreover, for each $uv\in E(F^{\bullet}[U,U'])$,
exactly one endpoint of $e(u,v)$ belongs to $S$.
This, together with $e(F^{\bullet}[U,U'])=2k$ and $w\in S$, implies that $|S|\geq 2k+1$,
which contradicts \eqref{align-9Y}.
Therefore, $H$ is $\mathcal{C}(F)$-free and $G=H$, as desired.

\medskip
\noindent{\textbf{Case 2.} $F^{\bullet}$ is the disjoint union of edges.}

Then, $F^{\bullet}\cong M_k$ for some $k\geq 1$.
By Claim \ref{claim-2.0}, we get $q=e(F^{\bullet})+1=k+1$.
Combining Claims \ref{claim-2.1} and \ref{claim-2.2},
we deduce that $G$ is a subgraph of $K_{k}\nabla E_{n-k}$.
To prove ${\rm EX}(n,\mathcal{M}(F^{o}))=\{K_{k}\nabla E_{n-k}\}$, it suffices to show that $K_{k}\nabla E_{n-k}$ is $\mathcal{C}(F)$-free.
Otherwise, $K_{k}\nabla E_{n-k}$ contains a member $J\in \mathcal{C}(F)$.
Then $\beta(J)\leq \beta(K_{k}\nabla E_{n-k})=k$.
However, this is impossible as $\beta(J)\geq \beta(F)=k+1$.

\medskip
\noindent{\textbf{Case 3.} All other cases for $F^{\bullet}$.}

Then $F^{\bullet}$ contains tree components, which implies that $e(F^{\bullet})\leq |F^{\bullet}|-1$.
By Claim \ref{claim-2.0}, we get $q=e(F^{\bullet})+1$.
Combining Claims \ref{claim-2.1} and \ref{claim-2.2},
we can see that $K_{e(F^{\bullet}),n-e(F^{\bullet})}\subseteq G$.
To prove ${\rm EX}(n,\mathcal{M}(F^{o}))=\{K_{e(F^{\bullet}),n-e(F^{\bullet})}\}$,
it suffices to show that $G[R]\cong E_{e(F^{\bullet})}$.
If every vertex of $A$ is of degree one in the graph $F^{\bullet}$, then $F^{\bullet}\cong M_k$ for some $k\geq 1$, a contradiction.
Thus $A$ contains a vertex, say $v_1$, such that $d_{F^{\bullet}}(v_1)\geq 2$.
If $G[R]$ is not empty,
then we can embed  $J_1$ into $G$ by letting $E(G[R])\cap E(J_1)=\{wv_1\}$, a contradiction.
Thus, $G[R]$ is empty, as desired.

This completes the proof of Lemma \ref{lem2.2}.
\end{proof}

\section{Some technical lemmas} \label{sec3}

The progressive induction, introduced by Simonovits \cite{Simonovits1968},
is a method used to prove statements that hold only for sufficiently large $n$.
It is similar to mathematical induction and the Euclidean algorithm,
combining elements of both in a certain sense.

\begin{lem}\emph{(\cite{Simonovits1968})}\label{lem3.2}
Let $\mathfrak{U}=\cup_{i=1}^{\infty}\mathfrak{U}_i$ be a set of given elements,
such that $\mathfrak{U}_i$ are disjoint finite subsets of $\mathfrak{U}$.
Let $B$ be a condition or property defined on $\mathfrak{U}$
(i.e., the elements of $\mathfrak{U}$ may satisfy or not satisfy B).
Let $\varphi(a)$ be a function defined on $\mathfrak{U}$ such that $\varphi(a)$ is a non-negative integer and
\begin{enumerate}[{\rm (i)}]\setlength{\itemsep}{0pt}
\item if $a$ satisfies $B$, then $\varphi(a)=0$;
\item there is an $N_0$ such that if $n>N_0$ and $a\in \mathfrak{U}_n$,
 then either $a$ satisfies $B$ or there exist an $n^*$ and an $a^*$ such that
$\frac{n}{2}<n^*<n$, $a^*\in \mathfrak{U}_{n^*}$, and $\varphi(a)<\varphi(a^*)$.
\end{enumerate}
Then there exists a positive integer $n_0$ such that for any $n > n_0$, every $a \in \mathfrak{U}_n$ satisfies $B$.

\end{lem}

The following lemmas will be crucial tools in the proof of our main result,
and we prove them here for later use.

\begin{lem}\label{lem3.3}
Let $H$ be a bipartite graph with color classes $A$ and $B$ satisfying the following:
\begin{enumerate}[{\rm (i)}]\setlength{\itemsep}{0pt}
\item $d_{B}(u)\geq \delta |B|$ for each $u\in A$ and some constant $0<\delta<1$,
where $d_{B}(u)=|N_H(u)\cap B|$;
\item $|A|\geq \frac{3\varphi}{\delta}$ for some integer $\varphi\geq 2$,
and $|B|$ is sufficiently large.
\end{enumerate}
Then there exist a $\varphi$-subset $A'\subseteq A$ and a $\lfloor(\frac{2\delta}{3})^{\varphi}|B|\rfloor$-subset $B'\subseteq B$ such that $H[A',B']\cong K_{|A'|,|B'|}$.
\end{lem}

\begin{proof}
Since $\varphi\geq 2$ and $0<\delta<1$, we get
$$\frac{3\varphi-k}{\lceil\frac{3\varphi}{\delta}\rceil-k}> \frac{3\varphi-\varphi}{\lceil\frac{3\varphi}{\delta}\rceil-\varphi}
>\frac{2\varphi}{\frac{3\varphi}{\delta}}=\frac{2\delta}{3}$$
for any $k\in \{0,1,\dots,\varphi-1\}$. Consequently,
\begin{align}\label{equ-5A}
\frac{\binom{3\varphi}{\varphi}}{\binom{\lceil\frac{3\varphi}{\delta}\rceil}{\varphi}}
=\frac{3\varphi}{\lceil\frac{3\varphi}{\delta}\rceil}\cdot \frac{3\varphi-1}{\lceil\frac{3\varphi}{\delta}\rceil-1} \cdots \frac{3\varphi-\varphi+1}{\lceil\frac{3\varphi}{\delta}\rceil-\varphi+1}
>\Big(\frac{2\delta}{3}\Big)^{\varphi}.
\end{align}

Let $A^*$ be a subset of $A$ with $|A^*|=\lceil\frac{3\varphi}{\delta}\rceil$ and
let $H^*=H[A^*,B]$.
By (i), we have
\begin{equation*}
e(H^*) = \sum_{u \in A^*} d_B(u) \geq |A^*| \cdot \delta |B| \geq 3\varphi |B|.
\end{equation*}
Let $a = \lfloor e(H^*) / |B| \rfloor$, which implies $a \geq 3\varphi$.
Among all sequences of non-negative integers $(d_1, \dots, d_{|B|})$ such that $d_1 \geq d_2 \geq \dots \geq d_{|B|}$ and $\sum_{j=1}^{|B|} d_j = e(H^*)$, we choose a sequence $(d_1^*, \dots, d_{|B|}^*)$ that minimizes the sum $\sum_{j=1}^{|B|} \binom{d_j}{\varphi}$.

We now claim that  $d_1^*-d_{|B|}^*\leq 1$.
Suppose, for the sake of contradiction, that $d_1^* - d_{|B|}^* \geq 2$,
implying $d_1^* - 1 > d_{|B|}^*$.
Since $d_1^*$ is the maximum value in a sequence with average $e(H^*)/|B|$,
we have $d_1^* \geq a \geq 3\varphi$, which yields $\binom{d_1^*-1}{\varphi-1} > \binom{d_{|B|}^*}{\varphi-1}$.
Using the identity $\binom{n}{k} = \binom{n-1}{k} + \binom{n-1}{k-1}$, we observe that
\begin{align*}
{d_1^*\choose \varphi}+{d_{|B|}^*\choose \varphi}
&={{d_1^*-1}\choose \varphi}+{{d_1^*-1}\choose \varphi-1}+{{d_{|B|}^*+1}\choose \varphi}
-{d_{|B|}^*\choose \varphi-1}\\
&>{{d_1^*-1}\choose \varphi}+{{d_{|B|}^*+1}\choose \varphi}.
\end{align*}
This implies that replacing $d_1^*$ and $d_{|B|}^*$ with $d_1^*-1$ and $d_{|B|}^*+1$ would result in a smaller sum of binomial coefficients, contradicting the minimality of $(d_1^*, \dots, d_{|B|}^*)$.
Thus, $d_1^* - d_{|B|}^* \leq 1$ must hold.

Since $d_1^*-d_{|B|}^*\leq 1$, it follows that $d_j^*\geq a\geq 3\varphi$ for all $j \in \{1, \dots, |B|\}$. Combining this with \eqref{equ-5A}, we have
\begin{align*}
\sum_{j=1}^{|B|}\binom{d_j^*}{\varphi} \geq |B|\binom{3\varphi}{\varphi} \geq |B|\left(\frac{2\delta}{3}\right)^{\varphi}\binom{|A^*|}{\varphi}.
\end{align*}
Note that $e(H^*) = \sum_{v\in B} d_{A^*}(v)$. By the construction of the sequence $(d_1^*, \dots, d_{|B|}^*)$ which minimizes the sum of binomial coefficients, we obtain
\begin{align}\label{equ-6A}
\sum_{v\in B}\binom{d_{A^*}(v)}{\varphi} \geq \sum_{j=1}^{|B|}\binom{d^*_j}{\varphi} \geq |B|\left(\frac{2\delta}{3}\right)^{\varphi}\binom{|A^*|}{\varphi}.
\end{align}
Let $\mathcal{A}$ be the family of all $\varphi$-subsets of $A^*$. For each $A^{**}\in \mathcal{A}$, let $B_{A^{**}} = \{v\in B : A^{**}\subseteq N_{A^*}(v)\}$. We consider the incidence set $\mathbb{A} = \{(v, A^{**}) : v\in B, A^{**}\in \mathcal{A}, A^{**}\subseteq N_{A^*}(v)\}$. By double counting $|\mathbb{A}|$ and applying \eqref{equ-6A}, we have
\begin{align*}
\sum_{A^{**}\in \mathcal{A}}|B_{A^{**}}| = |\mathbb{A}| = \sum_{v\in B}\binom{d_{A^*}(v)}{\varphi} \geq |B|\left(\frac{2\delta}{3}\right)^{\varphi}\binom{|A^*|}{\varphi}.
\end{align*}
Since $|\mathcal{A}| = \binom{|A^*|}{\varphi}$, the pigeonhole principle implies that there exists some $A' \in \mathcal{A}$ such that $|B_{A'}| \geq \lfloor (\frac{2\delta}{3})^{\varphi}|B| \rfloor$. Fix such an $A'$, and let $B' \subseteq B_{A'}$ be a subset of size $\lfloor (2\delta/3)^{\varphi}|B| \rfloor$. By the definition of $B_{A'}$, every vertex $v \in B'$ is adjacent to all vertices in $A'$. Consequently, $H[A', B']$ induces a complete bipartite graph $K_{\varphi, \lfloor(2\delta/3)^{\varphi}|B|\rfloor}$, which is the desired subgraph. This completes the proof.
\end{proof}

\begin{lem}\label{lem3.4}
Let $n$ be sufficiently large and $\mathcal{F}$ be a family of finite graphs such that $\min_{F\in \mathcal{F}}\chi(F)=r+1\geq 3$.
Let $G$ be an $n$-vertex $\mathcal{F}$-free graph obtained from $T_{n,r}$
by adding and deleting at most $\varepsilon n^2$ edges,
where
\begin{align}\label{equ-7A}
\phi=8^{r-1}\max_{F\in \mathcal{F}}|F|~~\text{and}~~\varepsilon^{\frac13}<\min\Big\{\frac{1}{4^{8^{r}\phi}\times 8\phi^4},\frac{1}{16\phi^6}\Big\}.
\end{align}
Define $V_1,\ldots, V_r$ such that they form a partition $V(G)=\bigcup_{i=1}^{r}V_i$
where $\sum_{1\leq i<j\leq r}e(V_i,V_j)$ is maximized.
Set $W=\bigcup_{i=1}^{r}W_i$,
where $W_i=\{w\in V_i~|~d_{V_i}(w)\geq 2\varepsilon^{\frac13} n\}$.
Then
\begin{enumerate}[{\rm (i)}]\setlength{\itemsep}{0pt}
\item $\sum_{i=1}^{r}e(V_i)\leq \varepsilon n^2$ and $\big||V_i|-\frac{n}{r}\big|\leq \varepsilon^{\frac13}n$ for each $i\in [r]$;
\item $|P|\leq \varepsilon^{\frac13} n$, where  $P=\{v\in V(G)~|~d_G(v)\leq \big(\frac{r-1}{r}-6\varepsilon^{\frac13}\big)n\}$;
\item $d_{V_{i}}(w)\geq(\frac{1}{2r}-\phi\varepsilon^{\frac13})n$ for any $i\in [r]$ and any $w\in \bigcup_{j\in [r]\setminus \{i\}}(W_{j}\setminus P)$;
\item $d_{V_{i}}(u)\geq\big(\frac{1}{r}-2\phi\varepsilon^{\frac13}\big)n$ for any $i\in [r]$ and any $u\in \bigcup_{j\in [r]\setminus \{i\}}(V_{j}\setminus (W_j\cup P))$;
\item $|W\setminus P|\leq {\varepsilon^{-\frac12}}$.
\end{enumerate}
\end{lem}

\begin{proof}
(i) Since $G$ is obtained from $T_{n,r}$ by adding and deleting at most $\varepsilon n^2$ edges,
we get
\begin{equation}\label{eq2.2}
e(G)\geq e(T_{n,r})-\varepsilon n^2\geq \frac{r-1}{2r}n^2-\frac{r}{8}-\varepsilon n^2
>\frac{r-1}{2r}n^2-2\varepsilon n^2,
\end{equation}
and there exists a partition $V(G)=\bigcup_{i=1}^{r} U_i$ such that
$\sum_{i=1}^{r}e(U_i)\leq \varepsilon n^2$
and $\big\lfloor\frac nr\big\rfloor\leq |U_i|
\leq\big\lceil\frac nr\big\rceil$ for each $i\in [r]$.
By the definitions of $V_1,\ldots, V_r$,
it follows that
\begin{align}\label{align-8B}
\sum_{i=1}^{r}e(V_i)\leq \sum_{i=1}^{r}e(U_i)\leq \varepsilon n^2.
\end{align}

For every $i\in [r]$, by the definition of $W_i$ we get
$$2e(V_i)=\sum_{w\in V_i}d_{V_i}(w)\geq
\sum_{w\in W_i}d_{V_i}(w)\geq |W_i|\cdot 2\varepsilon^{\frac13} n.
$$
Combining this with \eqref{align-8B} gives
$$
\varepsilon n^2\geq \sum_{i=1}^{r}e(V_i)\geq \sum_{i=1}^{r}|W_i|\varepsilon^{\frac13} n
  =|W|\varepsilon^{\frac13} n.
$$
This yields that $|W|\leq \varepsilon^{\frac23} n$.

Set $\alpha=\max\{\big||V_i|-\frac{n}{r}\big|~|~i\in [r]\}$. We may assume without loss of generality that $\alpha=||V_1|-\frac{n}{r}|$.
From the Cauchy-Schwarz inequality, we know that $(r-1)\sum_{i=2}^{r}|V_i|^2\geq (\sum_{i=2}^{r}|V_i|)^2$.
Then
   $$2\sum_{2\leq i<j\leq r}|V_i||V_j|=\Big(\sum_{i=2}^{r}|V_i|\Big)^2-\sum_{i=2}^{r}|V_i|^2
   \leq \frac{r-2}{r-1}(n-|V_1|)^2.$$
Since $\sum_{1\leq i<j\leq r}|V_i||V_j|=|V_1|(n-|V_1|)+\sum_{2\leq i<j\leq r}|V_i||V_j|$,
we have
\begin{eqnarray*}
\sum_{1\leq i<j\leq r}|V_i||V_j|
\leq |V_1|(n-|V_1|)+\frac{r-2}{2(r-1)}(n-|V_1|)^2
= -\frac{r}{2(r-1)}\alpha^2+\frac{r-1}{2r}n^2,
\end{eqnarray*}
where the last equality holds as $\alpha=||V_1|-\frac{n}{r}|$.
Consequently,
\begin{eqnarray*}
e(G)
\leq \sum_{1\leq i<j\leq r}|V_i||V_j|+\sum_{i=1}^{r}e(V_i)
= -\frac{r}{2(r-1)}\alpha^2+\frac{r-1}{2r}n^2+\varepsilon n^2.
\end{eqnarray*}
Combining with \eqref{eq2.2} gives $\frac{r}{2(r-1)}\alpha^2<3\varepsilon n^2$.
Thus $\alpha
<\sqrt{6\varepsilon} n<\varepsilon^{\frac13}n$
as $\varepsilon<16^{-3}$ by \eqref{equ-7A}.

(ii)
Suppose to the contrary that $|P|>\varepsilon^{\frac13} n$.
Then there exists a subset $P'\subseteq P$ with $|P'|=\lfloor\varepsilon^{\frac13} n\rfloor$.
Combining with \eqref{eq2.2}, we obtain
\begin{align}\label{equ-8A}
 e(G-P')&\geq  e(G)-\sum_{v\in P'}d_G(v)
        \geq  \Big(\frac{r-1}{2r}-2\varepsilon\Big)n^2-\varepsilon^{\frac13} n\Big(\frac{r-1}{r}-6\varepsilon^{\frac13}\Big)n\nonumber\\
  &=    \Big(\frac{r-1}{2r}-2\varepsilon-\frac{r-1}{r}\varepsilon^{\frac13}+6\varepsilon^{\frac23}\Big)n^2
  >\frac{r-1}{2r}(1-2\varepsilon^{\frac13}+6\varepsilon^{\frac23})n^2.
\end{align}
Define $n_0:=|G-P'|$.
Clearly, $n_0< (1-\varepsilon^{\frac13})n+1$.
On the other hand, since $G-P'$ is $\mathcal{F}$-free,
the Erd\H{o}s-Stone-Simonovits theorem implies that
 $$e(G-P')\leq \frac{r-1}{2r}n_0^2+o(n_0^2)
 <\frac{r-1}{2r}(1-2\varepsilon^{\frac13}+2\varepsilon^{\frac23})n^2,$$
which contradicts \eqref{equ-8A}.
Thus $|P|\leq \varepsilon^{\frac13} n$.

(iii) Choose an integer $i\in [r]$ and a vertex $w\in W_{j}\setminus P$ for some ${j}\in [r]\setminus \{i\}$.
Since $V(G)=\bigcup_{k=1}^{r}V_{k}$ is a partition
such that $\sum_{1\leq i'<j'\leq r}e(V_{i'},V_{j'})$ is maximized, we have
$d_{V_{{j}}}(w)\leq d_{V_{{i}}}(w)$.
By (i), $|V_{k}|\leq \big(\frac{1}{r}+\varepsilon^{\frac13}\big)n$ for any $k\in [r]$.
Combining these with $d_G(w)> \big(\frac{r-1}{r}-6\varepsilon^{\frac13}\big)n$ (as $w\notin P$) and \eqref{equ-7A},  we deduce that
\begin{align*}
  2d_{V_{i}}(w)
  &\geq d_{V_{i}}(w)+d_{V_{j}}(w)
  = d_G(w)-\sum_{k\in [r]\setminus\{i,j\}}d_{V_{k}}(w)\\
  &\geq \big(\frac{r-1}{r}-6\varepsilon^{\frac13}\big)n- (r-2)\Big(\frac{1}{r}+\varepsilon^{\frac13}\Big)n
    \geq \big(\frac{1}{r}-2\phi\varepsilon^{\frac13}\big)n.
\end{align*}
Thus, $d_{V_{i}}(w)\geq (\frac{1}{2r}-\phi\varepsilon^{\frac13})n$.

(iv)
Let $V_k^{\circ}=V_k\setminus (W_k\cup P)$ for each $k\in [r]$.
Choose an integer $i\in [r]$
and a vertex $u\in V_{j}^{\circ}$ for some $j\in [r]\setminus \{i\}$.
Since $u\notin P\cup W_j$, we get
$d_G(u)> \big(\frac{r-1}{r}-6\varepsilon^{\frac13}\big)n$
and $d_{V_{j}}(u)<2\varepsilon^{\frac13} n$.
Combining these with (i) and \eqref{equ-7A}, we obtain that
\begin{align*}
  d_{V_{i}}(u)
    &= d_G(u)-d_{V_{j}}(u)-\sum_{k\in [r]\setminus\{i,j\}}d_{V_{k}}(u)\\
    &\geq \Big(\frac{r-1}{r}-6\varepsilon^{\frac13}\Big)n-2\varepsilon^{\frac13} n- (r-2)\Big(\frac{1}{r}+\varepsilon^{\frac13}\Big)n
    \geq \Big(\frac{1}{r}-2\phi\varepsilon^{\frac13}\Big)n.
\end{align*}

(v) By symmetry, it suffices to prove that
$|W_1\setminus P|\leq \frac{1}{r}\varepsilon^{-\frac12}$.
Indeed, otherwise,
we get $$|W_1\setminus P|
>\frac{1}{r}\varepsilon^{-\frac12}
\geq \frac{1}{\phi}\times 4\phi^3  \varepsilon^{-\frac13}
> 3\times8^{r-1}\phi\varepsilon^{-\frac13}$$ by \eqref{equ-7A}.
For any $w\in W_{1}\setminus P$,  by $|W|\leq \varepsilon^{\frac23} n$, (i) and (ii) we get
\begin{align*}
d_{V_1^{\circ}}(w)
&\geq d_{V_1}(w)-|W|-|P|
\geq 2\varepsilon^{\frac13}n-\varepsilon^{\frac23}n-\varepsilon^{\frac13}n\\
&\geq \varepsilon^{\frac13}\Big(\frac{1}{r}+\varepsilon^{\frac13}\Big)n
\geq \varepsilon^{\frac13}|V_1^{\circ}|.
\end{align*}
Applying Lemma \ref{lem3.3} with $\delta=\varepsilon^{\frac13}$ and $\varphi=8^{r-1}\phi$,
there exist subsets $W^1\subseteq W_1\setminus P$ and $V_1^{\circ\circ}\subseteq V_1^{\circ}$ such that
$|W^1|=8^{r-1}\phi$, $|V_1^{\circ\circ}|=\phi$,
 and $H[W^1,V_1^{\circ\circ}]$ is a complete bipartite graph.
For any $w\in W^{1}$,
by (iii) we have $d_{V_{k}}(w)\geq(\frac{1}{2r}-\phi\varepsilon^{\frac13})n$
for any $k\in [r]\setminus \{1\}$.
Combining this with $|W|\leq \varepsilon^{\frac23} n$, (i), (ii) and \eqref{equ-7A}, we obtain that
\begin{align*}
d_{V_{k}^{\circ}}(w)
&>d_{V_{k}}(w)-|W|-| P|
>\Big(\frac{1}{2r}-\phi\varepsilon^{\frac13}\Big)n-\varepsilon^{\frac23}n-\varepsilon^{\frac13}n\\
&>\frac38\Big(\frac{1}{r}+\varepsilon^{\frac13}\Big)n>\frac38|V_k^{\circ}|.
\end{align*}
Applying Lemma \ref{lem3.3} with $\delta=\frac{3}{8}$ and $\varphi=8^{r-2}\phi$,
there exist subsets $W^2\subseteq W^1$ and $V^{\circ\circ}_2\subseteq V_2^{\circ}$ such that
$|W^2|=8^{r-2}\phi$, $|V^{\circ\circ}_2|=\lfloor\frac{|V_2^{\circ}|}{4^{|W^2|}}\rfloor$,
 and $H[W^2,V^{\circ\circ}_2]$ is a complete bipartite graph.
By repeating this process, we can obtain sequences of subsets $W^2,W^3,\dots,W^r$ and ${V}^{\circ\circ}_2,{V}^{\circ\circ}_3,\dots,{V}^{\circ\circ}_r$ such that
$W^k\subseteq W^{k-1}$, $V^{\circ\circ}_k\subseteq V_k^{\circ}$,
$|W^k|=8^{r-k}\phi$, $|{V}^{\circ\circ}_k|=\lfloor\frac{|V_k^{\circ}|}{4^{|W^k|}}\rfloor$,
 and $H[W^k,V^{\circ\circ}_k]$ is a complete bipartite graph for each $k\in [r]\setminus \{1\}$.

\begin{claim}\label{claim4.3C}
Let $i\in [r]\setminus \{1\}$ and $j\in [\phi^2]$.
Let $\{u_1,\dots,u_{j}\}\subseteq \bigcup_{k\in [r]\setminus \{i\}}V^{\circ\circ}_k$.
Then there exist at least $\phi$ vertices in $V^{\circ\circ}_{i}$
that are adjacent to $u_1,\dots,u_{j}$.
\end{claim}

\renewcommand\proofname{\bf Proof of Claim \ref{claim4.3C}}
\begin{proof}
For any $k\in [j]$,
we have $u_k\in V_{k_0}^{\circ\circ}$ for some $k_0\neq i$.
By (iv) and $V^{\circ\circ}_{k_0}\subseteq V_{k_0}^{\circ}$, we know that $d_{V_{i}}(u_k)\geq\big(\frac{1}{r}-2\phi\varepsilon^{\frac13}\big)n$.
This, together with (i), gives $|V_{i}\setminus N_{V_{i}}(u_k)|\leq 3\phi \varepsilon^{\frac13} n$.
It is clear that $V^{\circ\circ}_{i}\setminus N_{V^{\circ\circ}_{i}}(u_k)\subseteq V_{i}\setminus N_{V_{i}}(u_k)$. Hence,
\begin{align*}
\Big|\bigcup_{k=1}^{j} \big(V^{\circ\circ}_{i}\setminus N_{V^{\circ\circ}_{i}}(u_k)\big)\Big|
 \leq \sum_{k=1}^{j}|V^{\circ\circ}_{i}\setminus N_{V^{\circ\circ}_{i}}(u_k)|
 \leq  3\phi^3\varepsilon^{\frac13} n.
\end{align*}
By $|W|\leq \varepsilon^{\frac23} n$, (i) and (ii),
we get $|V_i^{\circ}|\geq |V_i|-|W|-|P|\geq \frac{n}{2\phi}$.
Combining this with \eqref{equ-7A} gives
$$|{V}^{\circ\circ}_i|
=\Big\lfloor\frac{|V_i^{\circ}|}{4^{|W^i|}}\Big\rfloor
\geq \Big\lfloor\frac{|V_i^{\circ}|}{4^{8^{r-1}\phi}}\Big\rfloor
\geq 4\phi^3\varepsilon^{\frac13} n.$$
Consequently,
\begin{align*}
\Big|\bigcap_{k=1}^{j}N_{V^{\circ\circ}_{i}}(u_k)\Big|
=  |V^{\circ\circ}_{i}|-\Big|\bigcup_{k=1}^{j} \big(V^{\circ\circ}_{i}\setminus N_{V_{i}}(u_k)\big)\Big|
\geq \phi^3\varepsilon^{\frac13} n.
\end{align*}
Then there exist at least $\phi$ vertices in $V^{\circ\circ}_{i}$ that are adjacent to $u_1,\dots,u_{j}$.
\end{proof}

Set $\widetilde{V}_1=V_1^{\circ\circ}$.
Clearly, $|\widetilde{V}_1|=|V_1^{\circ\circ}|=\phi$.
Recursively applying Claim \ref{claim4.3C},
there exists a sequence of $\phi$-subsets $\widetilde{V}_2,\dots,\widetilde{V}_r$
such that  for any $k\in [r]\setminus\{1\}$,
$\widetilde{V}_k\subseteq V^{\circ\circ}_k$ and each vertex in $\widetilde{V}_k$ is adjacent to all vertices in $\bigcup_{i=1}^{k-1}\widetilde{V}_i$.
Note that $|W^r|=\phi$ and $G[W^r,\widetilde{V}_k]$ is a complete bipartite graph for each $k\in [r]$.
These indicates that $G[W^r\cup (\bigcup_{k=1}^{r}\widetilde{V}_k)]\cong T_{(r+1)\phi,r+1}$.
Since $\phi>\max_{F\in \mathcal{F}}|F^{\bullet}|$ and $\min_{F\in \mathcal{F}}\chi(F)=r+1$,
it follows that $G[W^r\cup (\bigcup_{k=1}^{r}\widetilde{V}_k)]$ contains a member of $\mathcal{F}$,
a contradiction.
Thus (v) holds and the proof of Lemma \ref{lem3.4} is complete.
\end{proof}

\begin{lem}\label{lem3.5}
Let $n$ and $N$ be sufficiently large integers with $2N<n$, and
$F=K_1\nabla F^{\bullet}$ where each component of $F^{\bullet}$ is either a non-trivial tree or an even cycle.
Let $G$ be an $n$-vertex $F^{o}$-free graph with a partition of vertices into three parts $V(G)=R\cup V_1\cup V_2$
satisfying the following:
\begin{enumerate}[(a)]\setlength{\itemsep}{0pt}
 \item There exist $V_1'\subseteq V_1$ and $V_2'\subseteq V_2$ such that $G[V_1'\cup V_2']=G[V_1',V_2']=T_{2N,2}$;
\item $|R|=q(\mathcal{C}(F))-1$ and each vertex of $R$ is adjacent to each vertex of $T_{2N,2}$;
\item For each $i\in [2]$, each vertex of $V_i\setminus V_i'$ is adjacent to each vertex of $V_{3-i}'$.
\end{enumerate}
Then we have the following statements:
\begin{enumerate}[{\rm (i)}]\setlength{\itemsep}{0pt}
\item If each component of $F^{\bullet}$ is an even cycle, then $e(G)\leq e(H(n,\frac12{|F^{\bullet}|},2))$.
Moreover, the equality holds if and only if $G\cong H(n,\frac12{|F^{\bullet}|},2)$;
\item If each component of $F^{\bullet}$ is an edge,
then $e(G)\leq e(K_{e(F^{\bullet})}\nabla T_{n-e(F^{\bullet}),2})$.
Moreover, the equality holds if and only if $G\cong K_{e(F^{\bullet})}\nabla T_{n-e(F^{\bullet}),2}$;
\item For all other cases of $F^{\bullet}$, $e(G)\leq e(E_{e(F^{\bullet})}\nabla T_{n-e(F^{\bullet}),2})$.
Moreover, the equality holds if and only if $G\cong E_{e(F^{\bullet})}\nabla T_{n-e(F^{\bullet}),2}$.
\end{enumerate}
\end{lem}

\begin{proof}
We only prove the validity of (i).
The proofs of (ii) and (iii) are similar and hence omitted here.
Note that each component of $F^{\bullet}$ is an even cycle.
By Claim \ref{claim-2.0}, we get $|R|=q(\mathcal{C}(F))-1=|F^{\bullet}|-1$.
We first prove that $G[V_i]$ is empty for each $i\in [2]$.
Otherwise, by symmetry, we may assume that $G[V_1]$ contains an edge $uv$.
Note that $G[R,V_1']$ is a complete bipartite graph.
Recall that $q(J_4)=|F^{\bullet}|$ and $J_4$ contains exactly $|B|$ edge components.
It follows that $G[R\cup V_1'\cup \{u,v\}]$ contains a copy of $J_4$.
By the definition of $\mathcal{M}(F^{o})$ and $J_4\in \mathcal{M}(F^{o})$,
$G[(R\cup V_1'\cup \{u,v\})\cup V_2']$ contains a copy of $F^{o}$, a contradiction.
Thus, $G[V_i]$ is empty for each $i\in [2]$.

Suppose that $G[R]$ is $P_3$-free.
Then we can embed  $J_1$ into $G[R\cup V_1']$
by letting $E(G[R])\cap E(J_1)=\{wv_1,wv_2\}$, where $v_1,v_2$ belong to one color class of $F^{\bullet}$.
By the definition of $\mathcal{M}(F^{o})$ and $J_1\in \mathcal{M}(F^{o})$,
$G[(R\cup V_1')\cup V_2']$ contains a copy of $F^{o}$, a contradiction.
Thus, $G[R]$ is $P_3$-free.
It follows that $G[R]\subseteq M_{k-1}\cup K_1$,
and hence $G$ is a subgraph of
$(M_{k-1}\cup K_1)\nabla K_{|V_1|,|V_2|}$.
Consequently,
\begin{align}\label{equ-9A}
e(G)\leq e((M_{k-1}\cup K_1)\nabla K_{|V_1|,|V_2|})
\leq e((M_{k-1}\cup K_1)\nabla T_{n-2k+1,2})=h(n,k,2).
\end{align}
From Lemma \ref{lem2.2} we know that $H(n,k,1)$ is $\mathcal{M}(F^{o})$-free.
Then by the definition of $\mathcal{M}(F^{o})$,
it follows that $H(n,k,2)$ is $F^{o}$-free.
If $e(G)=h(n,k,2)$,
then by \eqref{equ-9A} we know that $K_{|V_1|,|V_2|}\cong T_{n-2k+1,2}$
and $G\cong H(n,k,2)$, as desired.
\end{proof}

\section{Proof of Theorem \ref{thm1.2A}}\label{sec5}

We use $\delta(G)$ to denote the minimum degree of $G$.
The following well-known lemma
is due to Erd\H{o}s \cite{Erdos1966} and Simonovits \cite{Simonovits1968}.

\begin{lem} \emph{(\cite{Erdos1966,Simonovits1968})}\label{lem4.1}
Let $F$ be a graph with $\chi(F)\geq 3$.
For any given $\varepsilon>0$,
if $n$ is sufficiently large, then
$\delta(G)>(1-\frac{1}{\chi(F)-1}-\varepsilon)n$ for any $G\in {\rm EX}(n,F)$.
\end{lem}

Building upon Lemma \ref{lem4.1}, we now establish the proof of Theorem \ref{thm1.2A}.
\renewcommand\proofname{\bf Proof of Theorem \ref{thm1.2A}}
\begin{proof}
Let $G$ be any graph in ${\rm EX}(n,F^{o})$.
Recall that $\ell_{u,v}$ is the length of the substituted odd cycle associated with the edge $uv\in E(F)$.
Suppose that $G$ contains a subgraph $F$.
We present the following claim.

\begin{claim}\label{claim4.1E}
Let $u_0v_0$ be any edge of $F$ and let $\ell_{u_0,v_0}=2k+1$.
For any $R\subseteq V(G)\setminus \{u_0,v_0\}$ with $|R|\leq |F^{o}|^2$,
$G-R$ contains a $(2k+1)$-cycle $C$ with $u_0v_0\in E(C)$.
\end{claim}

\renewcommand\proofname{\bf Proof of Claim \ref{claim4.1E}}
\begin{proof}
Assume $\varepsilon<\frac{1}{9}$.
By Lemma \ref{lem4.1},
we get $\delta(G)\geq (1-\frac{1}{3}-\frac{1}{9})n=\frac{5}{9}n$.
Since $\ell_{u_0,v_0}\geq 5$, we have $k\geq 2$.
We construct a path $P = u_0 u_1 \dots u_{2k-2}$ greedily such that for each $i \in [2k-2]$,
$u_i \in N_G(u_{i-1}) \setminus (R \cup \{u_0, \dots, u_{i-1}\} \cup \{v_0\}).$
Since $\delta(G)\geq \frac{5}{9}n$, we get
\begin{align*}
|N_G(u_{2k-2})\cap N_G(v_0)|\geq |N_G(u_{2k-2})|+|N_G(v_0)|-n
\geq \frac{n}{9}>|R\cup \{u_0,\dots,u_{2k-2}\}\cup \{v_0\}|.
\end{align*}
Then there exists a vertex $u_{2k-1}\in N_G(u_{2k-2})\cap N_G(v_0)$
such that $u_{2k-1}\notin R\cup \{u_0,\dots,u_{2k-2}\}\cup \{v_0\}$.
Consequently, $G-R$ contains a $(2k+1)$-cycle
$C:=u_0u_1\cdots u_{2k-1}v_0$, as desired.
\end{proof}

By recursively applying Claim \ref{claim4.1E} to all edges of $F$,
we see that $G$ contains a copy of $F^{o}$,
which leads to a contradiction.
Hence $G$ must be $F$-free.
It follows that ${\rm ex}(n,F^{o})=e(G)\leq {\rm ex}(n,F)$.
On the other hand, since $F$ is a subgraph of $F^{o}$,
we also get ${\rm ex}(n,F)\leq {\rm ex}(n,F^{o})$.
Therefore, ${\rm ex}(n,F)= {\rm ex}(n,F^{o})$.

Since $F$ is a subgraph of $F^{o}$ and ${\rm ex}(n,F)= {\rm ex}(n,F^{o})$,
we assert that ${\rm EX}(n,F)\subseteq {\rm EX}(n,F^{o})$.
Conversely, any $G \in \operatorname{EX}(n, F^o)$ is necessarily $F$-free. Given that $\operatorname{ex}(n, F) = \operatorname{ex}(n, F^o)$, we conclude that $\operatorname{EX}(n, F^o) \subseteq \operatorname{EX}(n, F)$. Therefore, $\operatorname{EX}(n, F) = \operatorname{EX}(n, F^o)$, which completes the proof of Theorem \ref{thm1.2A}.
\end{proof}

\section{Proof of Theorem \ref{thm1.1}}\label{sec4}

Throughout this section, assume that $n$ is sufficiently large and choose an arbitrary graph $G_n\in {\rm EX}(n,F^{o})$.
We focus on the proof of Theorem \ref{thm1.1}
for the case that each component of $F^{\bullet}$ is an even cycle.
The proofs for the remaining cases are analogous and are thus omitted.

Assume that $F^{\bullet}\cong \cup_{i=1}^{j}C_{2k_i}$.
Let $k=\sum_{i=1}^{j}k_i$,
and let $A=\{v_1,\dots,v_{k}\}$ and $B=\{u_1,\dots,u_{k}\}$ denote the two color classes of $F^{\bullet}$.
By Claim \ref{claim-2.0}, we have $q=e(F^{\bullet})=2k$.
By the famous Erd\H{o}s-Stone-Simonovits theorem,
$G_n$ can be obtained from $T_{n,2}$
by adding and deleting at most $\varepsilon n^2$ edges,
where
\begin{align}\label{equ-10A}
\phi=8|F^{o}|~~\text{and}~~
\varepsilon^{\frac13}<\min\Big\{\frac{1}{4^{64\phi}\times 8\phi^4},\frac{1}{16\phi^6}\Big\}.
\end{align}

Define $V_1,V_2$ such that they form a partition $V(G_n)=V_1\cup V_2$
where $e(V_1,V_2)$ is maximized.
Applying Lemma \ref{lem3.4} with $r=2$, we get the following statements:
\begin{enumerate}[(I)]\setlength{\itemsep}{0pt}
\item $e(V_1)+e(V_2)\leq \varepsilon n^2$ and $\big||V_i|-\frac{n}{2}\big|\leq \varepsilon^{\frac13}n$ for each $i\in [2]$;
\item Let $P=\{v\in V(G_n)~|~d_{G_n}(v)\leq \big(\frac{1}{2}-6\varepsilon^{\frac13}\big)n\}.$
Then $|P|\leq \varepsilon^{\frac13} n$;
\item Let $W=W_1\cup W_2$,
where $W_i=\{w\in V_i~|~d_{V_i}(w)\geq 2\varepsilon^{\frac13} n\}$.
Then $|W\setminus P|\leq \varepsilon^{-\frac12}$;
\item For any $i\in [2]$ and any $w\in W_{3-i}\setminus P$, we have $d_{V_{i}}(w)\geq(\frac{1}{4}-\phi\varepsilon^{\frac13})n$;
\item Let $V^{\circ}_{j}=V_{j}\setminus (W\cup P)$ for each $j\in [2]$.
For any $i\in [2]$ and any $u\in V_{3-i}^{\circ}$, we have $d_{V_{i}}(u)\geq\big(\frac{1}{2}-2\phi\varepsilon^{\frac13}\big)n$.
\end{enumerate}

By (II) and (III), we get
\begin{align}\label{align-10E}
|W\cup P|=|W\setminus P|+|P|\leq \varepsilon^{-\frac12}+\varepsilon^{\frac13} n.
\end{align}

In the following, we shall prove several claims for the extremal graph $G_n$.
\begin{claim}\label{claim4.1}
Let $i_0\in [2]$ and $j_0\in [\phi^2]$. Then the following hold:
\begin{enumerate}[(a)]\setlength{\itemsep}{0pt}
\item Let $u_0\in (R_1\cup R_2)\setminus P$ where $R_i=\{w\in V_i~|~d_{V_i}(w)\geq 4\phi^3\varepsilon^{\frac13} n\}$, and let $\{u_1,\dots,u_{j_0}\}\subseteq V_{3-i_0}^{\circ}$.
Then there exist at least $4\phi$ vertices in $V_{i_0}^{\circ}$
that are adjacent to $u_0,u_1,\dots,u_{j_0}$;
\item Let $u_0\in R_{3-i_0}\setminus P$ and $\{u_1,\dots,u_{j_0}\}\subseteq V_{3-i_0}\setminus (P\cup R_{3-i_0})$.
Then there exist at least $\phi$ vertices in $V_{i_0}^{\circ}$
that are adjacent to $u_0,u_1,\dots,u_{j_0}$.
\end{enumerate}
\end{claim}
\renewcommand\proofname{\bf Proof of Claim \ref{claim4.1}}
\begin{proof}
($a$)
For any $k\in [j_0]$, property (V) implies that $d_{V_{i_0}}(u_k)\geq\big(\frac{1}{2}-2\phi\varepsilon^{\frac13}\big)n$.
Combined with (I), this yields
$|V_{i_0}\setminus N_{V_{i_0}}(u_k)|\leq 3\phi\varepsilon^{\frac13} n$.
Consequently,
\begin{align}\label{equ-11A}
\big|\bigcup_{k=1}^{j_0} \big(V_{i_0}\setminus N_{V_{i_0}}(u_k)\big)\big|
 \leq \sum_{k=1}^{j_0}|V_{i_0}\setminus N_{V_{i_0}}(u_k)|
 \leq  3\phi^3\varepsilon^{\frac13} n.
\end{align}
Since $u_0\in (R_1\cup R_2)\setminus P$,
property (IV) implies $$d_{V_{i_0}}(u_0)
\geq\min\Big\{\Big(\frac{1}{4}-\phi\varepsilon^{\frac13}\Big)n,4\phi^3\varepsilon^{\frac13} n\Big\}
=4\phi^3\varepsilon^{\frac13} n.$$
Combining this with \eqref{align-10E} and \eqref{equ-11A}, we obtain
\begin{align*}
\Big|N_{V_{i_0}}(u_0)\cap\Big(\bigcap_{k=1}^{j_0}N_{V_{i_0}}(u_k)\Big)\Big|
\geq |N_{V_{i_0}}(u_0)|-\Big|\bigcup_{k=1}^{j_0} \big(V_{i_0}\setminus N_{V_{i_0}}(u_k)\big)\Big|
\geq \phi^3\varepsilon^{\frac13} n
>|P\cup W|+4\phi,
\end{align*}
where the last inequality holds as $\varepsilon^{\frac13}<\frac{1}{16\phi^6}$.
Thus, there exist at least $4\phi$ vertices in $V_{i_0}^{\circ}$ that are adjacent to $u_0,u_1,\dots,u_{j_0}$, as desired.

($b$) For any $k\in [j_0]$, since $u_k\in V_{3-i_0}\setminus (P\cup R_{3-i_0})$,
we get $d_G(u_k)> \big(\frac{1}{2}-6\varepsilon^{\frac13}\big)n$ and $d_{V_{3-i_0}}(u_k)\leq 4\phi^3\varepsilon^{\frac13} n$.
Thus, $$d_{V_{i_0}}(u_k)\geq d_G(u_k)-d_{V_{3-i_0}}(u_k)\geq \big(\frac{1}{2}-6\varepsilon^{\frac13}-4\phi^3\varepsilon^{\frac13}\big)n.$$
This, together with (I), gives that
$|V_{i_0}\setminus N_{V_{i_0}}(u_k)|\leq 5\phi^3\varepsilon^{\frac13} n$.
Then,
\begin{align*}
\Big|\bigcup_{k=1}^{j_0} \big(V_{i_0}\setminus N_{V_{i_0}}(u_k)\big)\Big|
 \leq \sum_{k=1}^{j_0}|V_{i_0}\setminus N_{V_{i_0}}(u_k)|
 \leq  5\phi^5\varepsilon^{\frac13} n.
\end{align*}
Since $u_0\in R_{3-i_0}\setminus P\subseteq W_{3-i_0}\setminus P$, we have $d_{V_{i_0}}(u_0)\geq \big(\frac{1}{4}-\phi\varepsilon^{\frac13}\big)n$ by (IV).
Combining these with \eqref{equ-10A} and \eqref{align-10E}, we get that
\begin{align*}
\Big|N_{V_{i_0}}(u_0)\cap\Big(\bigcap_{k=1}^{j_0}N_{V_{i_0}}(u_k)\Big)\Big|
\geq |N_{V_{i_0}}(u_0)|-\Big|\bigcup_{k=1}^{j_0} \big(V_{i_0}\setminus N_{V_{i_0}}(u_k)\big)\Big|
>|P\cup W|+\phi.
\end{align*}
Thus, there exist at least $\phi$ vertices in $V_{i_0}^{\circ}$
that are adjacent to $u_1,\dots,u_{j_0}$, as desired.
\end{proof}

\begin{claim}\label{claim4.2}
$P=\varnothing$.
\end{claim}

\renewcommand\proofname{\bf Proof of Claim \ref{claim4.2}}
\begin{proof}
Suppose to the contrary that there exists a vertex $u_0\in P$.
By the definition of $P$, we have $d_{G_n}(u_0)\leq \big(\frac{1}{2}-6\varepsilon^{\frac13}\big)n.$
Let $G'_n$ be the graph obtained from $G_n$ by deleting all edges incident to $u_0$
and adding all possible edges between $V_1^{\circ}\setminus \{u_0\}$ and $u_0$.
By (I) and  \eqref{align-10E}, we get
$$e(G'_n)-e(G_n)\geq |V_1^{\circ}\setminus \{u_0\}|-d_G(u_0)\geq |V_1|-|P\cup W|-1-d_G(u_0)>0.$$
By the choice of $G_n$,
we can see that $G'_n$ contains a subgraph $H'$ isomorphic to $F^{o}$.
From the construction of $H'$, it follows that $u_0\in V(H')$.
Assume that $N_{H'}(u_0)=\{u_1,u_2,\dots,u_c\}$.
Then $c\leq |H'|-1<\phi$ and $u_1,\dots,u_c\in V_1^{\circ}$.
By Claim \ref{claim4.1} ($a$),
there exists a vertex $u'\in V_{2}^{\circ}\setminus V(H')$ adjacent to $u_1,u_2,\dots,u_c$.
This implies that $G_n[(V(H')\setminus \{u_0\})\cup\{u'\}]$ contains a copy of $H'$,
which is impossible since $G_n$ is $F^{o}$-free.
Thus, $P=\varnothing$.
\end{proof}

\begin{claim}\label{claim4.3}
For each $i\in [2]$, we have $\nu(G_n[\overline{V}_i])\leq \phi$,
where $\overline{V}_i=V_i\setminus R_i$.
\end{claim}

\renewcommand\proofname{\bf Proof of Claim \ref{claim4.3}}
\begin{proof}
By symmetry,
it suffices to prove that $\nu(G_n[\overline{V}_1])\leq \phi$.
Indeed, otherwise, $G_n[\overline{V}_1]$ contains $\phi$ independent edges, say $w_1w_2,\dots,w_{2\phi-1}w_{2\phi}$.
Set $V_1'=\{w_1,\dots,w_{2\phi}\}$.

We first prove that $R_1=\varnothing$.
Suppose for contradiction that $R_1\neq \varnothing$, and let $w_0\in R_1$.
By the definition of $R_1$,
we have $d_{V_1}(w_0)\geq 4\phi^3\varepsilon^{\frac13} n$.
Since $P=\varnothing$ and $R_1\subseteq W$,
property (III) implies that $|R_1|\leq |W\setminus P|\leq \varepsilon^{-\frac12}$.
Consequently,
$d_{\overline{V}_1}(w_0)\geq d_{V_1}(w_0)-|R_1|\geq 3\phi^3\varepsilon^{\frac13} n$.
This guarantees the existence of a $k$-subset $V_1''\subseteq N_{\overline{V}_1}(w_0)\setminus V_1'$.
Since $V_1'\cup V_1''\subseteq \overline{V}_1$ and $P=\varnothing$,
Claim \ref{claim4.1} ($b$) ensures that
there exists a $\phi$-subset $V_2'\subseteq V_2^{\circ}$ such that each vertex in $V_2'$ is adjacent to all vertices in $V_1'''$, where $V_1'''=V_1'\cup V_1''\cup \{w_0\}$.
It is easy to find that $G_n[V_1''']$ contains a copy of $J_3$.
By the definition of $\mathcal{M}(F^{o})$,
$G_n[V_1'''\cup V_2']$ contains a copy of $F^{o}$, a contradiction.
Therefore, we must have $R_1 = \varnothing$, which in turn implies $\overline{V}_1 = V_1$.

Suppose that $\Delta(G_n[V_1])> k$.
Then there exist a vertex $w_0\in V_1$ and a $k$-subset $V_1''\subseteq N_{V_1}(w_0)$.
Note that $\phi=8|F^{o}|>16k$.
Then there exist at least $3k$ edges in $\{w_1w_2,\dots,w_{2\phi-1}w_{2\phi}\}$
 that are incident to no vertex of $\{w_0\}\cup V_1''$.
Consequently, $G_n[V_1''']$ contains a copy of $J_2$,
where $V_1'''=(\{w_0\}\cup V_1'')\cup V_1'$.
Note that $V_1'''\subseteq V_1=\overline{V}_1$ and $P=\varnothing$.
Then by Claim \ref{claim4.1} ($b$),
there exists a $\phi$-subset $V_2'\subseteq V_2^{\circ}$ such that
each vertex in $V_2'$ is adjacent to all vertices in $V_1'''$.
From the definition of $\mathcal{M}(F^{o})$ we know that
$G_n[V_1'''\cup V_2']$ contains a copy of $F^{o}$, a contradiction.
Thus, $\Delta(G_n[V_1])\leq k$.
It follows that $W_1=\varnothing$ and $V_1=V_1^{\circ}$.

Since $H(n,k,2)$ is $F^{o}$-free,
we have $e(G_n)\geq h(n,k,2)$.
Then,
\begin{align}\label{equ-12A}
e(V_1)+e(V_2)=e(G_n)-e(V_1,V_2)\geq h(n,k,2)-e(T_{n,2})= \frac{2k-1}{2}n-O(1).
\end{align}
Let $w$ be a vertex of $G_n[V_1]\cup G_n[V_2]$ with maximum degree $\Delta$.
Then $\Delta\geq \lceil2\frac{\frac{2k-1}{2}n-O(1)}{n}\rceil=2k-1$.
This, together with $\Delta(G_n[V_1])\leq k$, implies that $w\in V_2$.
On the other hand,
$e(V_1)\leq \frac{1}{2}\Delta(G_n[V_1])|V_1|\leq \frac{k}{2}(\frac{1}{2}+\varepsilon^{\frac13})n$.
Combining this with \eqref{equ-12A}, we obtain
\begin{align}\label{equ-13A}
e(V_2)\geq \frac{k-1}{2}n.
\end{align}

In what follows, we split the proof into the following two cases.

\medskip
\noindent{\textbf{Case 1.} $\Delta< 4\phi^3\varepsilon^{\frac13} n$.}

Since $\Delta(G_n[V_2])=d_{V_2}(w)\geq 2k-1$,
there is a $k$-subset $V_2''\subseteq N_{V_2}(w)$.
Suppose that $\nu(G_n[V_2])\geq \phi$.
That is, $G_n[V_2]$ contains $\phi$ independent edges, say $w_1'w_2',\dots,w_{2\phi-1}'w_{2\phi}'$.
Set $V_2'=\{w_1',\dots,w_{2\phi}'\}$.
Note that $\phi=8|F^{o}|>16k$.
Then there exist at least $3k$ edges in $\{w_1'w_2',\dots,w_{2\phi-1}'w_{2\phi}'\}$
 that are incident to no vertex in $\{w\}\cup V_2''$.
 Clearly, $J_2\cong S_{k+1}\cup M_{3k}$.
One can further observe that $G_n[V_2''']$ contains a copy of $J_2$,
where $V_2'''=(\{w\}\cup V_2'')\cup V_2'$.
Since $\Delta< 4\phi^3\varepsilon^{\frac13} n$,
it is clear that $R_2=\varnothing$ and $V_2'''\subseteq \overline{V}_2$.
By Claim \ref{claim4.1} ($b$) and $P=\varnothing$,
there exists a $\phi$-subset $V_1'\subseteq V_1^{\circ}$ such that
each vertex in $V_1'$ is adjacent to all vertices in $V_2'''$.
From the definition of $\mathcal{M}(F^{o})$ we know that
$G_n[V_2'''\cup V_1']$ contains a copy of $F^{o}$, a contradiction.
Thus, $\nu(G_n[V_2])< \phi$.
Combining this with Lemma \ref{lem3.1} and \eqref{equ-10A} gives that
$$e(V_2)\leq \nu(G_n[V_2])(\Delta(G_n[V_2])+1)
\leq(\phi-1)(4\phi^3\varepsilon^{\frac13} n+1)
<\frac{k-1}{2}n,$$
which contradicts \eqref{equ-13A}.

\medskip
\noindent{\textbf{Case 2.} $\Delta\geq 4\phi^3\varepsilon^{\frac13} n$.}

Clearly, $w\in R_2$.
Since $W_1=\varnothing$ and $P=\varnothing$,
we have $V_1=V_1^{\circ}$, which implies $V_1'\subseteq V_1^{\circ}$.
Combining these with Claim \ref{claim4.1} ($a$) and $P=\varnothing$,
there exists a $\phi$-subset $V_2''\subseteq V_2^{\circ}$ such that each vertex in $V_2''$ is adjacent to all vertices in $V_1'\cup \{w\}$.
Again by Claim \ref{claim4.1} ($a$) and $P=\varnothing$,
there exists a $k$-subset $V_1''\subseteq V_1^{\circ}\setminus V_1'$ such that each vertex in $V_1''$ is adjacent to all vertices in $V_2''\cup \{w\}$.
It is clear that $G_n[V_1'\cup (V_1''\cup\{w\})]$ contains a copy of $J_2$.
By the definition of $\mathcal{M}(F^{o})$ and $J_2\in \mathcal{M}(F^{o})$,
we can deduce that $G_n[(V_1'\cup V_1''\cup \{w\})\cup V_2'']$ contains a copy of $F^{o}$, a contradiction.

This completes the proof of Claim \ref{claim4.3}.
\end{proof}

\begin{claim}\label{claim4.4}
There exist subsets $V_1^*\subseteq \overline{V}_1$ and $V_2^*\subseteq \overline{V}_2$ such that
$G_n[V_1^*,V_2^*]\cong T_{2\lfloor n^{0.8}\rfloor,2}$, where $V_i^*$ is an independent set of $G_n[\overline{V}_i]$ for each $i\in [2]$.
\end{claim}

\renewcommand\proofname{\bf Proof of Claim \ref{claim4.4}}
\begin{proof}
For $i\in [2]$,
let $\{u_{i,1}u_{i,2},\dots,u_{i,2\phi_i-1}u_{i,2\phi_i}\}$
be a maximum matching in $G_n[\overline{V}_i]$.
By Claim \ref{claim4.3}, we get  $\phi_i\leq \phi$.
Set $$V_i^{**}=\{u_{i,1},u_{i,2},\dots,u_{i,2\phi_i}\}
\cup \{v~|~v\in \bigcup\limits_{j=1}^{2\phi_i}N_{\overline{V}_i}(u_{i,j})\}.$$
Clearly, $\overline{V}_i\setminus V_i^{**}$ is an independent set of $G_n[\overline{V}_i]$.
For each $i\in [2]$ and each $j\in [2\phi_i]$,
since $u_{i,j}\in \overline{V}_i$,
we get $d_{\overline{V}_i}(u_{i,j})<4\phi^3\varepsilon^{\frac13} n$.
Consequently,
$|V_i^{**}|\leq (4\phi^3\varepsilon^{\frac13} n+1)\times 2\phi=9\phi^4\varepsilon^{\frac13} n$.


Suppose that $d_{V_{i}}(u_0)< |V_{i}|-2\varepsilon^{-\frac12}$
for some vertex $u_0\in \overline{V}_{3-i}$.
Let $G'_n$ be the graph obtained from $G_n$ by deleting all edges incident to $u_0$
and adding all possible edges between $u_0$ and $V_{i}^{\circ}$.
By Claim \ref{claim4.2} and (III), we get $|W|\leq \varepsilon^{-\frac12}$.
Thus
$$e(G'_n)-e(G_n)\geq (|V_{i}|-\varepsilon^{-\frac12})-d_{V_{i}}(u_0)>0.$$
By the choice of $G_n$,
it holds that $G'_n$ contains a subgraph $H'$ isomorphic to $F^{o}$.
From the construction of $H'$, we know that $u_0\in V(H')$.
Assume that $N_{H'}(u_0)=\{u_1,u_2,\dots,u_c\}$.
By Claim \ref{claim4.1} ($a$) and $P=\varnothing$,
there exists a vertex $u'\in V^{\circ}_{3-i}\setminus V(H')$ adjacent to $u_1,u_2,\dots,u_c$.
This implies that $G_n[(V(H')\setminus \{u_0\})\cup\{u'\}]$ contains a copy of $H'$,
which is impossible since $G_n$ is $F^{o}$-free.
Thus, $d_{V_{i}}(u)\geq |V_{i}|-2\varepsilon^{-\frac12}$
for each $u\in \overline{V}_{3-i}$.

Select a $\lfloor n^{0.8}\rfloor$-subset $V_1^*\subseteq \overline{V}_1\setminus V_1^{**}$.
Then,
$\big|\bigcup_{u\in V_1^*}(V_2\setminus N_{V_2}(u))\big|
\leq 2\varepsilon^{-\frac12}\lfloor n^{0.8}\rfloor$.
By (III) and Claim \ref{claim4.2}, we get $|R_2|\leq |W|\leq \varepsilon^{-\frac12}$.
Combining these with $|V_2^{**}|\leq 9\phi^4\varepsilon^{\frac13} n$, (i) and \eqref{equ-10A}, we obtain that
$$\Big|\bigcap\limits_{u\in V_1^*}N_{V_2}(u)\Big|
\geq |V_2|-2\varepsilon^{-\frac12}\lfloor n^{0.8}\rfloor
\geq \lfloor n^{0.8}\rfloor+|V_2^{**}|+|R_2|.$$
Then there exists a $\lfloor n^{0.8}\rfloor$-subset
$V_2^*\subseteq \overline{V}_2\setminus V_{2}^{**}$ such that each vertex in $V_2^*$ is adjacent to all vertices in $V_1^*$.
Therefore, $V_1^*$ and $V_2^*$ are the desired subsets.
\end{proof}

With Claims \ref{claim4.1}-\ref{claim4.4} established,
we are ready to complete the proof of Theorem \ref{thm1.1}
for the case that each component of $F^{\bullet}$ is an even cycle.

\renewcommand\proofname{\bf Proof of Theorem \ref{thm1.1}}
\begin{proof}
By Lemma \ref{lem3.5} (i), $H(n,k,2)$ is $F^{o}$-free.
Let $G_n$ be an arbitrary graph in ${\rm EX}(n,F^{o})$.
Then
\begin{align}\label{equ-14A}
e(G_n)={\rm ex}(n,F^{o})\geq e(H(n,k,2))=h(n,k,2).
\end{align}
Let $\mathfrak{U}_n={\rm EX}(n,F^{o})$
and $B$ be the property defined on $\mathfrak{U}=\cup_{n=1}^{\infty}\mathfrak{U}_n$ stating that $G_n\cong H(n,k,2)$.
Define $\varphi(G_n)=e(G_n)-h(n,k,2)$.
By \eqref{equ-14A}, $\varphi(G_n)$ is non-negative,
and it is 0 whenever $G_n$ satisfies property $B$.
Hence, Condition (i) of Lemma \ref{lem3.2} is satisfied.
To complete the proof of Theorem \ref{thm1.1},
it suffices to verify  Condition (ii) of Lemma \ref{lem3.2}.
Specifically, we need to show that $G_n\cong H(n,k,2)$ or there exists an integer $n^*\in (\frac{n}{2},n)$
such that $\varphi(G_{n^*})>\varphi(G_{n})$.
Once Condition (ii) is established,
Lemma \ref{lem3.2} guarantees an $n_0$ such that $G_n\cong H(n,k,2)$ for all $n>n_0$,
thereby completing the proof of Theorem \ref{thm1.1}.

In what follows, we will find a subgraph of $G_n$ satisfying the conditions of Lemma \ref{lem3.5}.
Let $L_0=G_n[V_1^*\cup V_2^*]$,
where $V_1^*$ and $V_2^*$ are defined in Claim \ref{claim4.4}.
For $i\geq 1$,
we will define vertices $w_i\in R$ and graphs $L_i$ recursively.
If there exists a vertex $w\in R\setminus \{w_1,\dots, w_{i-1}\}$
which has at least $ \varepsilon^{2i} n^{0.8}$ neighbors in each color class of $L_{i-1}$,
then let $w_i=w$ and define $L_i=T_{2\lfloor \varepsilon^{2i}n^{0.8}\rfloor,2}$ as a subgraph of $L_{i-1}$ in the neighbors of $w_i$.
Recall that $q=q(\mathcal{M}(F^{o}))$.
If $i\geq q$ , then by the definition of $q$,
$G_n[\{x_1,\dots,x_q\}\cup V(L_q)]$ contains a complete $3$-partite subgraph $K_{q,\lfloor \varepsilon^{2q}n^{0.8}\rfloor,\lfloor \varepsilon^{2q}n^{0.8}\rfloor}$, and hence contains a copy of $F^{o}$,
which leads to a contradiction.
Thus, $i\leq q-1$.

Now suppose the above process ends at $L_{\ell}$ for some $\ell\leq q-1$.
Set $R^*=\{w_1,\dots,w_\ell\}$.
Clearly, each vertex in $R^*$ is adjacent to all vertices of $L_\ell$.
Let $B_1^{\ell}$ and $B_2^{\ell}$ be the color classes of $L_\ell$,
where $B_1^{\ell}\subseteq V_1^*$ and $B_2^{\ell}\subseteq V_2^*$.
Then, $|B_i^\ell|=\lfloor \varepsilon^{2\ell} n^{0.8}\rfloor$ for each $i\in [2]$.
Partition the vertices $x\in V(G_n)\setminus (V(L_\ell)\cup R^*)$ into the following vertex sets:
\begin{itemize}\setlength{\itemsep}{0pt}
\item If there exists an $i\in [2]$ such that $x$ is adjacent to less than $\varepsilon^{2\ell+2}n^{0.8}$
vertices of $B_i^{\ell}$ and is adjacent to at least $(1-\varepsilon)\varepsilon^{2\ell} n^{0.8}$ vertices of $B_{3-i}^{\ell}$, then let $x\in D_i$;
\item If there exists an $i\in [2]$ such that $x$ is adjacent to less than $\varepsilon^{2\ell+2}n^{0.8}$
vertices of $B_i^{\ell}$ and is adjacent to less than $(1-\varepsilon)\varepsilon^{2\ell} n^{0.8}$ vertices of $B_{3-i}^{\ell}$, then let $x\in D$.
\end{itemize}
Obviously, $D_1\cup D_2\cup D$ is a partition of $V(G_n)\setminus (V(L_\ell)\cup R^*)$.
Select a maximal set of independent edges in $G_n[B_i^\ell,D_i]$,
say $x_1y_1,\dots,x_ty_t$, where $x_j\in B_i^\ell$, $y_j\in D_i$ for each $j\in [t]$.
By Claim \ref{claim4.4} and $B_i^{\ell}\subseteq V_i^*$, $B_i^{\ell}$ is an independent set of $G_n[\overline{V}_i]$.
It follows that $D_i\subseteq (\overline{V}_i\setminus B_i^\ell)\cup (R\setminus R^*)$
and $y_j\in R\setminus R^*$ for each $j\in [t]$.
By Claim \ref{claim4.2} and (III), $t\leq |R|\leq |W|\leq \varepsilon^{-\frac12}$.
Let $B_i^*=\{v\in B_i^\ell~|~N_{D_i}(v)\neq \varnothing\}$.
Since $x_1y_1,\dots,x_ty_t$ is a maximal set of independent edges,
we get $B_i^*= \cup_{j=1}^{t}N_{B_i^\ell}(y_j)$.
By the definition of $D_i$,
we know that each vertex $y_j$ is adjacent to less than $\varepsilon^{2\ell+2}n^{0.8}$ of $B_i^\ell$,
it holds that $|B_i^*|\leq \varepsilon^{2\ell+\frac32}n^{0.8}$.
Select a subset $B_i^{**}\subseteq B_i^\ell$
such that $B_i^{*}\subseteq B_i^{**}$ and $|B_i^{**}|=\lfloor\varepsilon^{2\ell+\frac32}n^{0.8}\rfloor$.
Set  $B_i'=B_i^\ell\setminus B_i^{**}$ and $D_i'=D_i\cup B_i^{**}$.
By Claim \ref{claim4.4} and $B_i'\subseteq B_i^\ell\subseteq V_i^*$,
we know that $B_i'$ is an independent set of $G_n[\overline{V}_i]$.
Thus, there are no edges between $B_i'$ and $D_i'$.
Let $L_{\ell}'=G_n[B_1'\cup B_2']$ and $\widetilde{G}=G_n-V(L_{\ell}')$.
Clearly, $R^*\cup D_1'\cup D_2'\cup D$ is a partition of $V(\widetilde{G})$, and we observe that:
\begin{itemize}\setlength{\itemsep}{0pt}
\item Each vertex of $R^*$ is adjacent to all vertices of $L_{\ell}'$;
\item Each vertex of $D_i'$  is adjacent to no vertex of $B_i'$ and is adjacent to at least $(1-\varepsilon-\varepsilon^{\frac32})\varepsilon^{2\ell}n^{0.8}$ vertices of $B_{3-i}'$;
\item Each vertex of $D$ is adjacent to less than $\varepsilon^{2\ell+2}n^{0.8}$ vertices of $B_i'$ and
is adjacent to less than $(1-\varepsilon)\varepsilon^{2\ell}n^{0.8}$ vertices of $B_{3-i}'$ for some $i\in [2]$.
\end{itemize}

Denote by $e_L$ the number of edges between $V(\widetilde{G})$ and $V(L_{\ell}')$.
Clearly, we have
\begin{align}\label{equ-15A}
e(G_n)=e(L_{\ell}')+e_L+e(\widetilde{G}).
\end{align}
By Claim \ref{claim4.4} and $B_i'\subseteq B_i^\ell\subseteq V_i^*$,
we obtain $L_{\ell}'=G_n[B_1'\cup B_2']\cong T_{2n',2}$, where $n'=|B_i'|$.
Select an induced copy of $L_\ell'$ in the graph $H(n,k,2)$.
Let $\widetilde{H}=H(n,k,2)-V(L_\ell')$ and $e_H$ be the number of edges between $V(\widetilde{H})$ and $V(L_{\ell}')$.
Then, $\widetilde{H}\cong H(n-2n',k,2)$ and
\begin{align}\label{equ-16A}
e(H(n,k,2))=e(L_{\ell}')+e_H+e(\widetilde{H}).
\end{align}
Since $\widetilde{G}$ is $F^{o}$-free, we have $e(\widetilde{G})\leq e(G_{n-2n'})$,
where $G_{n-2n'}\in {\rm EX}(n-2n',F^{o})$.
By \eqref{equ-15A} and \eqref{equ-16A}, we obtain that
\begin{align*}
\varphi(G_n)
&=e(G_n)-e(H(n,k,2))=(e_L-e_H)+(e(\widetilde{G})-e(\widetilde{H}))\\
& \leq (e_L-e_H)+(e(G_{n-2n'})-e(H(n-2n',k,2)))
 \leq (e_L-e_H)+\varphi(G_{n-2n'}).
\end{align*}
By \eqref{equ-10A}, we get $\varepsilon^{\frac12}+\varepsilon<1$,
and hence $1-\varepsilon^{\frac32}>1-\varepsilon+\varepsilon^{2}$. Thus,
\begin{align*}
n'=|B_i^\ell\setminus B_i^{**}|
=\lfloor \varepsilon^{2\ell} n^{0.8}\rfloor-\lfloor \varepsilon^{2\ell+\frac32}n^{0.8}\rfloor
>(1-\varepsilon)\varepsilon^{2\ell}n^{0.8}+\varepsilon^{2\ell+2}n^{0.8}.
\end{align*}
Consequently,
\begin{align*}
e_L
&\leq \ell\cdot 2n'+(n-\ell-2n'-|D|)n'+|D|((1-\varepsilon)\varepsilon^{2\ell}n^{0.8}+\varepsilon^{2\ell+2}n^{0.8})\\
&\leq  \ell\cdot 2n'+(n-\ell-2n')n'+|D|((1-\varepsilon)\varepsilon^{2\ell}n^{0.8}+\varepsilon^{2\ell+2}n^{0.8}-n')\\
&\leq (n+\ell-2n')n'.
\end{align*}
Furthermore, if $e_L=(n+\ell-2n')n'$, then $|D|=0$,
and each vertex of $D_i'$ is adjacent to all vertices of $B'_{3-i}$ for any $i\in [2]$.
On the other hand, $$e_H=(q-1)2n'+(n-q+1-2n')n'=(n+q-1-2n')n'.$$
Since $\ell\leq q-1$, we have
$e_L-e_H\leq (\ell-q+1)n'\leq 0$.
Set $n^*=n-2n'$.
If $e_L-e_H<0$, then $\varphi(G_n)<\varphi(G_{n^*})$.
Since $n'=|B_i'|=\lfloor \varepsilon^{2\ell} n^{0.8}\rfloor-
\lfloor \varepsilon^{2\ell+\frac32}n^{0.8}\rfloor=\Theta(n^{0.8})$, it follows that $n^*=n-2n'>\frac{n}{2}$ and we are done.
It remains to deal with the case that $e_L=e_H$.
This implies that $|D|=0$, $\ell=q-1$,
and each vertex of $D_i'$ is adjacent to all vertices of $B_{3-i}'$  for any $i\in [2]$.
Applying Lemma \ref{lem3.5} to $G_n$,
we have  $e(G_n)\leq h(n,k,2)$ and  the equality holds if and only if $G_n\cong H(n,k,2)$.
This, together with \eqref{equ-14A}, gives that $G_n\cong H(n,k,2)$.
This  completes the proof of Theorem \ref{thm1.1}.
\end{proof}


\begin{thebibliography}{99}
\setlength{\itemsep}{0pt}


\bibitem {Abbott-1972}
H.L. Abbott, D. Hanson, N. Sauer,
Intersection theorems for systems of sets,
\emph{J. Combin. Theory Ser. A}
\textbf{12} (1972), 381--389.

\bibitem {Chen-2003}
G.T. Chen, R.J. Gould, F. Pfender, B. Wei,
Extremal graphs for intersecting cliques,
\emph{J. Combin. Theory Ser. B.}
\textbf{89} (2) (2003), 159--171.

\bibitem {CH-1976}
V. Chv\'{a}tal, D. Hanson, Degrees and matchings,
\emph{J. Combin. Theory Ser. B}
\textbf{20} (2) (1976), 128--138.


\bibitem {E1962}
P. Erd\H{o}s, \"{U}ber ein Extremalproblem in der Graphentheorie (German), \emph{Arch. Math. (Basel)} \textbf{13} (1962), 122--127.

\bibitem {Erdos1966}
P. Erd\H{o}s,
On some new inequalities concerning extremal properties of graphs,
Theory of Graphs (Proc. Colloq., Tihany, 1966), pp. 77--81, Academic Press, New York-London, 1968.




\bibitem {EFGG1995}
P. Erd\H{o}s, Z. F\"{u}redi, R.J. Gould, D.S. Gunderson,
Extremal graphs for intersecting triangles,
\emph{J. Combin. Theory Ser. B} \textbf{64} (1) (1995), 89--100.


\bibitem {Erdos-1966}
P. Erd\H{o}s, M. Simonovits,
A limit theorem in graph theory,
\emph{Studia Sci. Math. Hungar.}
\textbf{1} (1966), 51--57.

\bibitem {ES1946}
P. Erd\H{o}s, A.H. Stone,
On the structure of linear graphs,
\emph{Bull. Amer. Math.}
\textbf{52} (1946), 1087--1091.




\bibitem {Hou2016}
X.M. Hou, Y. Qiu, B.Y. Liu,
Extremal graph for intersecting odd cycles,
\emph{Electron. J. Combin.}
\textbf{23} (2) (2016), Paper 2.29, 9 pp.

\bibitem {Hou2018}
X.M. Hou, Y. Qiu, B.Y. Liu,
Tur\'{a}n number and decomposition number of intersecting odd cycles,
\emph{Discrete Math.}
\textbf{341} (1) (2018), 126--137.


\bibitem {LIU}
H. Liu,
Extremal graphs for blow-ups of cycles and trees,
\emph{Electron. J. combin.}
\textbf{20} (1) (2013), Paper 65, 16pp.

\bibitem {Moon1968}
J.W. Moon,
On independent complete subgraphs in a graph,
\emph{Canad. J. Math.}
\textbf{20} (1968), 95--102.


\bibitem {Peng-2024}
X. Peng, M.J. Xia,
Tur\'{a}n number of the odd-ballooning of complete bipartite graphs,
\emph{J. Graph Theory} \textbf{107} (1) (2024), 181--199.

\bibitem {Simonovits1968}
M. Simonovits,
A method for solving extremal problems in graph theory, stability problems,
\emph{in: Theory of Graphs (Proc. Colloq., Tihany, 1966)}, \emph{Academic Press, New York,}
1968, pp. 279--319.

\bibitem {Simonovits-1983}
M. Simonovits,
Extremal graph problems and graph products,
\emph{Studies in Pure Mathematics,
Akad\'{e}miai Kiad\'{o} and Birkh\"{a}user Verlag}, (1983), 669--680.


\bibitem {Simonovits-1974}
M. Simonovits,
Extremal graph problems with symmetrical extremal graphs,
additional chromatic conditions, \emph{Discrete Math.}
\textbf{7} (1974), 349--376.

\bibitem{Simonovits-2019}
M. Simonovits, E. Szemer\'edi,
Embedding graphs into larger graphs: results, methods, and problems, In: {\em Building Bridges II: Mathematics of L\'aszl\'o Lov\'asz}, (2019), 445-592.


\bibitem {Song2024}
J.L. Song, C.H. Lu, L.-T. Yuan,
On the Tur\'{a}n number of edge blow-ups of cliques,
\emph{SIAM J. Discrete Math. }
\textbf{38} (3) (2024), 2429--2446.


\bibitem {Turan}
 P. Tur\'{a}n,
 On an extremal problem in graph theory,
 \emph{Mat. Fiz. Lapok} \textbf{48} (1941), 436--452.

\bibitem {Wang2021}
A.Y. Wang, X.M. Hou, B.Y. Liu, Y. Ma,
The Tur\'{a}n number for the edge blow-up of trees,
\emph{Discrete Math.} \textbf{344} (12) (2021), Paper No. 112627, 12 pp.

\bibitem {Yuan2018}
L.-T. Yuan,
Extremal graphs for the $k$-flower,
\emph{J. Graph Theory}
\textbf{89} (1) (2018),  26--39.

\bibitem {Y2022}
L.-T. Yuan,
Extremal graphs for edge blow-up of graphs,
\emph{J. Combin. Theory Ser. B} \textbf{152} (2022), 379--398.

\bibitem {Zhu2020}
H. Zhu, L.Y. Kang, E.F. Shan,
Extremal graphs for odd-ballooning of paths and cycles,
\emph{Graphs Combin.}
\textbf{36} (3) (2020), 755--766.

\bibitem {Zhu2023}
X.T. Zhu, Y.J. Chen,
Tur\'{a}n number for odd-ballooning of trees,
\emph{J. Graph Theory} \textbf{104} (2) (2023), 261--274.

\end{thebibliography}
\end{document}